\documentclass[aap,MSNbibl,seceqn,noautosecdot,dvips]{arximspdf}

\doi{10.1214/13-AAP980} 
\volume{24}
\issue{6}
\pubyear{2014}
\firstpage{2371}
\lastpage{2413}

\makeatletter
\newcommand{\widebar}{\overline}

\newcommand{\rrVert}{\Vert}
\newcommand{\llVert}{\Vert}
\newtheorem{theorem}{Theorem}[section]%
\newtheorem{corollary}[theorem]{Corollary}%
\newtheorem{proposition}[theorem]{Proposition}%
\newtheorem{lemma}[theorem]{Lemma}%
\newcommand{\rL}{\mathrm{L}}
\newcommand{\R}{\mathbb{R}}
\newcommand{\IND}{\mathbf{1}}
\makeatother

\begin{document}
\begin{frontmatter}

\title{First-order global asymptotics for confined particles with
singular pair repulsion}
\runtitle{First-order global asymptotics for confined particles}

\begin{aug}
\author[a]{\fnms{Djalil} \snm{Chafa\"{i}}\ead[label=e1]{djalil@chafai.net}\ead[label=u1,url]{http://djalil.chafai.net/}},
\author[b]{\fnms{Nathael} \snm{Gozlan}\ead[label=e3]{Nathael.Gozlan@u-pem.fr}}
\and
\author[b]{\fnms{Pierre-Andr\'e} \snm{Zitt}\corref{}\ead[label=e2]{Pierre-Andre.Zitt@u-pem.fr}}
\runauthor{D. Chafa\"{i}, N. Gozlan and P.-A. Zitt}
\affiliation{Universit\'e Paris-Est Marne-la-Vall\'ee}
\address[a]{D. Chafa\"{i}\\
Laboratoire d'Analyse\\
\qquad et de Math\'ematiques Appliqu\'ees\\
CNRS UMR8050\\
Universit\'e Paris-Est Marne-la-Vall\'ee\\
5 bd Descartes, F-77454, Champs-sur-Marne\\
France\\
Current address:\\
Universit\'e Paris-Dauphine\\
CEREMADE, IUF, PSL\\
France\\
\printead{e1}\\
\printead{u1}} 
\address[b]{N. Gozlan\\
P.-A. Zitt\\
Laboratoire d'Analyse\\
\qquad et de Math\'ematiques Appliqu\'ees\\
CNRS UMR8050\\
Universit\'e Paris-Est Marne-la-Vall\'ee\\
5 bd Descartes, F-77454, Champs-sur-Marne\\
France\\
\printead{e3}\\
\phantom{E-mail:}\ \printead*{e2}}
\end{aug}

\received{\smonth{6} \syear{2013}}
\revised{\smonth{10} \syear{2013}}

%
\begin{abstract}
We study a physical system of $N$ interacting particles in $\mathbb{R}^d$,
$d\geq1$, subject to pair repulsion and confined by an external field. We
establish a large deviations principle for their empirical
distribution as
$N$ tends to infinity. In the case of Riesz interaction, including Coulomb
interaction in arbitrary dimension $d>2$, the rate function is strictly
convex and admits a unique minimum, the equilibrium measure, characterized
via its potential. It follows that almost surely, the empirical distribution
of the particles tends to this equilibrium measure as $N$ tends to infinity.
In the more specific case of Coulomb interaction in dimension $d>2$, and
when the external field is a convex or increasing function of the radius,
then the equilibrium measure is supported in a ring. With a quadratic
external field, the equilibrium measure is uniform on a ball.
\end{abstract}

%

\end{frontmatter}

\section{\texorpdfstring{Introduction.}{Introduction}}

We study in this work a physical system of $N$ particles at positions
$x_1,\ldots,x_N\in\mathbb{R}^d$, $d\geq1$, with identical
``charge'' $q_N:=1/N$,
subject to a confining potential $V\dvtx\mathbb{R}^d\to\mathbb{R}$
coming from an
external field
and acting on each particle, and to an interaction potential
$W\dvtx\mathbb{R}^d\times\mathbb{R}^d\to(-\infty,+\infty]$
acting on each pair of
particles. The
function $W$ is finite outside the diagonal and symmetric: for all
$x,y\in\mathbb{R}^d$ with $x\neq y$, we\vspace*{1pt} have $W(x,y)=W(y,x)<\infty$.
The energy
$H_N(x_1,\ldots,x_N)$ of the configuration $(x_1,\ldots,x_N)\in
(\mathbb{R}^d)^N$
takes the form
%
%
\begin{eqnarray}\label{eqdefHn}
&& H_N(x_1,\ldots,x_N)\nonumber
\\
&&\qquad:= \sum_{i=1}^Nq_NV(x_i)+\sum_{i<j}q_N^2W(x_i,x_j)
\nonumber\\[-8pt]\\[-8pt]
&&\qquad =\frac{1}{N}\sum_{i=1}^NV(x_i)%
+\frac{1}{N^2}\sum_{i<j}W(x_i,x_j)\nonumber
\\
&&\qquad = \int V(x) \,d\mu_N(x)+\frac{1}{2}\int\!\!\!\int
_{\neq} W(x,y) \,d\mu_N(x) \,d\mu_N(y),\nonumber
\end{eqnarray}
where $\mu_N:=\frac{1}{N}\sum_{i=1}^N\delta_{x_i}$ is the empirical
measure of
the particles, and where the subscript ``$\neq$'' indicates that the double
integral is off-diagonal. The energy $H_N\dvtx(\mathbb{R}^d)^N\to
\mathbb{R}\cup\{
+\infty\}
$ is a
quadratic form functional in the variable~$\mu_N$.

From now on, and unless otherwise stated, we denote by ${{\vert
\cdot\vert}}$ the
Euclidean norm of $\mathbb{R}^d$, and we make the following additional
assumptions:
\begin{longlist}[(H3)]
\item[(H1)] The function $W\dvtx\mathbb{R}^d\times\mathbb
{R}^d\to
(-\infty,+\infty]$
is continuous on $\mathbb{R}^d\times\mathbb{R}^d$, symmetric, takes
finite values on
$\mathbb{R}^d\times\mathbb{R}^d \setminus\{(x,x); x\in\mathbb
{R}^d\}$ and satisfies the
following integrability condition: for all compact subset $K\subset
\mathbb{R}^d$,
the function
\[
z\in\mathbb{R}^d\mapsto\sup\bigl\{W(x,y); {{\vert x-y \vert}}\geq{{
\vert z \vert}}, x,y\in K \bigr\}
\]
is locally Lebesgue-integrable on $\R^d$.
\item[(H2)] The function $V\dvtx\mathbb{R}^d\to\R$ is
continuous and such
that\break  $\lim_{|x| \to+ \infty} V(x)=+\infty$ and
\[
\int_{\mathbb{R}^d}%
\exp{{ \bigl(-V(x) \bigr)}} \,dx<\infty.
\]
\item[(H3)] There exist constants $c\in\mathbb{R}$ and
$\varepsilon_o \in
(0,1)$ such that for every \mbox{$x,y\in\mathbb{R}^d$},
\[
W(x,y)\geq c-\varepsilon_o \bigl(V(x)+V(y) \bigr).
\]
(This must be understood as ``$V$ dominates $W$ at infinity.'')
\end{longlist}
Let ${(\beta_N)}_{N}$ be a sequence of positive real numbers such that
$\beta_N\to+\infty$ as \mbox{$N\to\infty$}. Under \textup{(H2)--(H3)}, there
exists an
integer $N_0$ depending on $\varepsilon_o$ such that for any $N\geq
N_0$, we have
\[
Z_N:=\int_{\mathbb{R}^d}\cdots\int_{\mathbb{R}^d}
\exp{{ \bigl(-\beta_NH_N(x_1,
\ldots,x_N) \bigr)}} \,dx_1\cdots dx_N<\infty,
\]
so that we can define the Boltzmann--Gibbs probability measure $P_N$ on
$(\mathbb{R}^d)^N$ by
%
%
\begin{equation}
\label{eqPN} dP_N(x_1,\ldots,x_N):=
\frac{\exp{{ (-\beta_N H_N(x_1,\ldots,x_N) )}}}{Z_N} \,dx_1\cdots dx_N.
\end{equation}
The law $P_N$ is the equilibrium distribution of a system of $N$ interacting
Brownian particles in $\mathbb{R}^d$, at inverse temperature $\beta
_N$, with equal
individual ``charge'' $1/N$, subject to a confining potential $V$
acting on
each particle, and to an interaction potential $W$ acting on each pair of
particles; see Section~\ref{sssSDE}. Note that for \mbox{$\beta_N=N^2$},
the quantity
$\beta_NH_N$ can also be interpreted as the distribution of a system
of $N$
particles living in $\mathbb{R}^d$, with unit ``charge,'' subject to a
confining
potential $NV$ acting on each particle, and to an interaction potential $W$
acting on each pair of particles.

Our work is motivated by the following physical control problem: given the
(internal) interaction potential $W$, for instance, a Coulomb
potential, a
target probability measure $\mu_\star$ on $\mathbb{R}^d$, for
instance, the uniform
law on the unit ball, and a cooling scheme $\beta_N\to+\infty$, for instance,
$\beta_N=N^2$, can we tune the (external) confinement potential $V$ (associated
to an external confinement field) such that $\mu_N\to\mu_\star$ as
$N\to\infty$? In this direction, we provide some partial answers in
Theorems~\ref{thldp}~and~\ref{thriesz}, Corollaries~\ref{cocoulomb}~and~\ref{coprescription} below. We also discuss several possible
extensions and related problems in Section~\ref{sscomments-extensions-related}.

Let $\mathcal{M}_1(\mathbb{R}^d)$ be the set of probability measures
on $\mathbb{R}
^d$. The
mean-field symmetries
of the model suggest to study, under the exchangeable measure $P_N$, the
behavior as $N\to\infty$ of the empirical measure $\mu_N$, which is
a random
variable on $\mathcal{M}_1(\mathbb{R}^d)$. With this asymptotic
analysis in
mind, we
introduce the functional $I\dvtx\mathcal{M}_1(\mathbb{R}^d)\to
(-\infty,+\infty
]$ given by
\[
I(\mu):= \frac{1}{2}\int\!\!\!\int{{ \bigl(V(x)+V(y)+W(x,y) \bigr
)}} \,d
\mu(x)\,d\mu(y).
\]
Assumptions \textup{(H2)}--\textup{(H3)} imply that the function
under the
integral is bounded from below, so that the integral defining $I$ makes sense
in $\mathbb{R}\cup\{+\infty\}=(-\infty,+\infty]$. If it is finite, then
$\int
V\,d\mu$
and $\int\!\!\int W\,d\mu^2$ both exist (see Lem\-ma~\ref{lemgoodrate}),
so that
\[
I(\mu) = \int V \,d\mu+ \frac{1}{2} \int\!\!\!\int W \,d\mu^2.
\]
The energy $H_N$ defined by (\ref{eqdefHn}) is ``almost'' given by
$I(\mu_N)$, where the infinite terms on the diagonal are forgotten.

\subsection{\texorpdfstring{Large deviations principle.}{Large deviations principle}}

Theorem~\ref{thldp} below is our first main result. It is of topological
nature, inspired from the available results for logarithmic Coulomb
gases in
random matrix theory \cite{MR1465640,MR1660943,MR1606719,MR2926763}.
We equip
$\mathcal{M}_1(\mathbb{R}^d)$ with the weak topology,\vspace*{1pt} defined by duality
with bounded
continuous functions. For any set $A\subset\mathcal{M}_1(\mathbb
{R}^d)$ we
denote by
$\operatorname{int}(A)$, $\operatorname{clo}(A)$ the interior and
closure of $A$ with respect to this
topology. This topology can be metrized by the Fortet--Mourier distance
defined by (see \cite{MR0061325,MR1619170})
%
%
\begin{equation}
\label{eqFortetMourier} {d_\mathrm{FM}}(\mu,\nu):=%
\sup
_{\max(|f|_\infty,|f|_{\mathrm{Lip}})\leq1}{{ \biggl\{\int f \,
d\mu
-\int
f \,d\nu\biggr\}}},
\end{equation}
where $|f|_\infty:=\sup|f|$ and $|f|_{\mathrm{Lip}}:=\sup_{x\neq
y}|f(x)-f(y)|/|x-y|$.

To formulate the large deviations result we need to introduce the following
additional technical assumption:
\begin{longlist}[(4)]
\item[(H4)] For all $\nu\in\mathcal{M}_1(\mathbb{R}^d)$
such that
$I(\nu)<+\infty$, there is a sequence $(\nu_n)_{n\in\mathbb{N}}$
of probability
measures, absolutely continuous with respect to Lebesgue, such that
$\nu_n$
converges weakly to $\nu$ and $I(\nu_n) \to I(\nu)$, when $n\to
\infty$.
\end{longlist}
It turns out that assumption~\textup{(H4)} is satisfied for a large
class of
potentials $V,W$, and several examples are given in
Proposition~\ref{propassumptionH4} and Theorem~\ref{thriesz}.

Throughout the paper, if $(a_N)_{N}$ and $(b_N)_{N}$ are nonnegative
sequences, the notation $a_N \gg b_N$ means that $a_N=b_Nc_N$, for some
$c_N$ that goes to $+\infty$ when $N\to\infty$.
%
%
\begin{theorem}[(Large deviations principle)]\label{thldp}
Suppose that
\[
\beta_N\gg N\log(N).
\]
If \textup{(H1)--(H3)} are satisfied, then:
\begin{longlist}[(4)]
\item[(1)] $I$ has compact level sets (and is thus lower semi-continuous) and
$\inf_{\mathcal{M}_1(\mathbb{R}^d)}I>-\infty$.
\item[(2)]  Under $(P_N)_N$, the sequence ${(\mu_N)}_{N}$ of random elements
of $\mathcal{M}_1(\mathbb{R}^d)$ equipped with the weak topology has
the following
asymptotic properties. For every Borel subset $A$ of $\mathcal
{M}_1(\mathbb{R}^d)$,
\begin{eqnarray*}
\limsup_{N\to\infty}\frac{\log Z_NP_N(\mu_N\in A)}{\beta_N} &\leq&
-\inf
_{\mu\in\operatorname{clo}(A)}I(\mu)
\end{eqnarray*}
and
\begin{eqnarray*}
\liminf_{N\to\infty}\frac{\log Z_NP_N(\mu_N\in A)}{\beta_N} &\geq&
-\inf
\bigl\{I(\mu); \mu
\in\operatorname{int}(A), \mu\ll\mathrm{Lebesgue} \bigr\}.
\end{eqnarray*}
\item[(3)]  Under the additional assumption $\textup{(H4)}$, the full Large
Deviation Principle (LDP) at speed $\beta_N$ holds with the rate function
\[
I_\star:=I-\inf_{\mathcal{M}_1(\mathbb{R}^d)}I.
\]
More precisely, for all Borel set $A \subset\mathcal{M}_1(\mathbb{R}^d)$,
\begin{eqnarray*}
-\inf_{\mu\in\operatorname{int}(A)} I_\star(\mu) &\leq&\liminf
_{N\to\infty} \frac{\log P_N(\mu_N \in A)}{\beta_N}
\\
&\leq&\limsup_{N\to\infty} \frac{\log P_N(\mu_N \in A)}{\beta_N}
\leq
-\inf
_{\mu\in\operatorname{clo}(A)} I_\star(\mu).
\end{eqnarray*}
In particular, by taking $A=\mathcal{M}_1(\mathbb{R}^d)$, we get
\[
\lim_{N\to\infty}\frac{\log Z_N}{\beta_N} =\inf_{\mathcal
{M}_1(\mathbb{R}
^d)}I_\star.
\]
\item[(4)]  Let $I_{\min}:=\{\mu\in\mathcal{M}_1\dvtx I_\star(\mu)=0\}
\neq
\varnothing
$. If
$\textup{(H4)}$ is satisfied and if ${(\mu_N)}_{N}$ are constructed
on the same probability space, and if $d$ stands for the Fortet--Mourier
distance (\ref{eqFortetMourier}), then we have, almost surely,
\[
\lim_{N\to\infty}{d_\mathrm{FM}}(\mu_N,I_{\min})=0.
\]
\end{longlist}
\end{theorem}

A careful reading of the proof of Theorem~\ref{thldp} indicates that if
$I_{\min}=\{\mu_\star\}$ is a singleton, and if \textup{(H4)} holds for
$\nu=\mu_\star$, then $\mu_N\to\mu_\star$ almost surely as $N\to
\infty$.

\subsection{\texorpdfstring{Case $\beta_N=N$ and link with Sanov theorem.}{Case beta N = N and link with Sanov theorem}}

If we set $W=0$, then the particles become i.i.d., and $P_N$ becomes a product
measure $\eta_N^{\otimes N}$ where $\eta_N\propto e^{-(\beta_N/N)V}$,
where the
symbol ``$\propto$'' means ``proportional to.'' When $\beta_N=N$, then
$\eta_N\propto e^{-V}$ does not depend on $N$, and we may denote it
$\eta$. To
provide perspective, recall that the classical Sanov theorem \cite
{MR2571413}, Theorem~6.2.10, for i.i.d. sequences means in our settings
that if $W=0$
and $\beta_N=N$, then ${(\mu_N)}_N$ satisfies to a large deviations
principle on
$\mathcal{M}_1(\mathbb{R}^d)$ at speed $N$ and with good rate function
\[
\mu\mapsto K(\mu|\eta):= %
\cases{ \displaystyle\int f\log(f) \,d\eta, &
\quad if $\mu\ll\eta$, with $\displaystyle f:=\frac{d\mu}{d\eta}$;
\cr
+\infty, &
\quad otherwise} %
\]
(Kullback--Leibler relative entropy or free energy). This large deviations
principle corresponds to the convergence
$\lim_{N\to\infty}{d_\mathrm{FM}}(\mu_N,\eta)=0$. Note that if
$\mu$ is absolutely
continuous with respect to Lebesgue measure with density function $g$, then
$K(\mu|\eta)$ can be decomposed in two terms,
\[
K(\mu|\eta) = \int V \,d\mu-H(\mu)+\log Z_V,
\]
where $Z_V:=\int_{\mathbb{R}^d} e^{-V(x)} \,dx$ and where $H(\mu)$
is the
Boltzmann--Shannon ``continuous'' entropy $H(\mu):=
-\int g(x)\log(g(x)) \,dx$; therefore at the speed $\beta_N = N$,
the energy
factor $\int V \,d\mu$ and the Boltzmann--Shannon entropy factor
$H(\mu
)$ both
appear in the rate function. In contrast, note that Theorem~\ref{thldp}
requires a higher inverse temperature $\beta_N\gg N\log(N)$. If we set
$W=0$ in
Theorem~\ref{thldp}, then $P_N$ becomes a product measure, the
particles are
i.i.d. though their common law depends on $N$, the function $\mu
\mapsto
I_*(\mu)=\int V \,d\mu-\inf V$ is affine, its minimizers $I_{\min}$ over
$\mathcal{M}_1(\mathbb{R}^d)$ coincide with
\[
\mathcal{M}_V:= \bigl\{\mu\in\mathcal{M}_1 \bigl(
\mathbb{R}^d \bigr)\dvtx\operatorname{supp}(\mu)\subset\arg\inf
V \bigr
\}
\]
and Theorem~\ref{thldp} boils down to a sort of Laplace principle, which
corresponds to the convergence $\lim_{N\to\infty}{d_\mathrm
{FM}}(\mu
_N,\mathcal{M}_V)=0$.
It is worthwhile to notice that the main difficulty in Theorem~\ref{thldp}
lies in the fact that $W$ can be infinite on the diagonal (short scale
repulsion). If $W$ is continuous and bounded on $\mathbb{R}^d\times
\mathbb{R}^d$,
then one
may deduce the large deviations principle for ${(\mu_N)}_{N}$ from the case
$W=0$ by using the Laplace--Varadhan lemma \cite{MR2571413},
Theorem~4.3.1; see
also \cite{MR1465640}, Corollary~5.1. To complete the picture, let us mention
that if $\beta_N=N$ and if $W$ is bounded and continuous, then the
Laplace--Varadhan lemma and the Sanov theorem would yield to the conclusion
that $(\mu_N)_N$ verifies a large deviations principle on $\mathcal
{M}_1(\mathbb{R}
^d)$ at
speed $N$ with rate function $R-\inf_{\mathcal{M}_1(\mathbb{R}^d)}R$
where the
functional $R$
is defined by
\begin{eqnarray*}
R(\mu) &:=& K(\mu|\eta) + \frac{1}{2}\int\!\!\!\int W(x,y) \,d\mu
(x)\,d\mu(y)
\\
&=& -H(\mu) + I(\mu)+\log Z_V;
\end{eqnarray*}
once more, the Boltzmann--Shannon entropy factor $H(\mu)$ reappears at
this rate. For an alternative point of view, we refer to
\cite{MR704588}, \cite{MR1145596}, Theorem~2.1, \cite{MR1362165,MR1678526,Kie93} and
\cite{KieSpo99}.

\subsection{\texorpdfstring{Equilibrium measure.}{Equilibrium measure}}

Our second-main result, expressed in Theorem~\ref{thriesz} and
Corollary~\ref{cocoulomb} below is of differential nature. It is based
on an
instance of the general Gauss problem in potential theory \cite
{frostman,MR0350027,MR2075132,MR2073873}. It concerns special choices
of $V$ and $W$
for which $I_\star$ achieves its minimum $0$ for a unique and explicit
$\mu_\star\in\mathcal{M}_1(\mathbb{R}^d)$. Recall that the \emph
{Coulomb} interactions
correspond to the choice $W(x,y)=k_\Delta(x-y)$ where $k_\Delta$ is the
\emph{Coulomb kernel} (opposite in sign to the Newton kernel) defined on
$\mathbb{R}^d$, $d\geq1$, by
%
%
\begin{equation}
\label{eqcoul} k_\Delta(x):= %
\cases{ -{{\vert x \vert
}}, &\quad if $d=1$,
\vspace*{4pt}\cr
\displaystyle\log\frac{1}{{{\vert x \vert}}}, &\quad if $d=2$,
\vspace*{4pt}\cr
\displaystyle\frac{1}{{{\vert x \vert}}^{d-2}},&\quad if $d\geq
3$.} %
\end{equation}
This is, up to a multiplicative constant, the fundamental
solution\footnote{There are no boundary conditions here, and thus the term
``Green function'' is not appropriate.
} of the Laplace equation. In other
words, denoting $\Delta:={\partial}_{x_1}^2+\cdots+{\partial
}_{x_d}^2$ the Laplacian, we have,
in a weak sense, in the space of Schwartz--Sobolev distributions
$\mathcal{D}
'(\mathbb{R}^d)$,
%
%
\begin{equation}
\label{eqfs-coul} -c\Delta k_\Delta=\delta_0\qquad\mbox{with }
c:= %
\cases{ \displaystyle\frac{1}{2}, &\quad if $d=1$,
\vspace*{4pt}\cr
\displaystyle\frac{1}{2\pi}, &\quad if $d=2$,
\vspace*{4pt}\cr
\displaystyle
\frac{1}{d(d-2)\omega_d}, &\quad if $d\geq3$,} %
\end{equation}
where $\omega_d:=\frac{\pi^{d/2}}{\Gamma(1+d/2)}$ is the volume of
the unit ball
of $\mathbb{R}^d$. Our notation is motivated by the fact that $-\Delta
$ is a
nonnegative operator.
The case of Coulomb interactions in dimension $d=2$ is known as ``logarithmic
potential with external field'' and is widely studied in the
literature; see
\cite{MR1746976,MR1485778,MR2760897,MR2926763}. To focus on novelty,
we will
not study the Coulomb kernel for $d\leq2$. We refer to
\cite{MR0148411,MR0147214,MR0129874,brascamp-lieb,MR597033,sandier-serfaty-1d}
and references therein for the Coulomb case in dimension $d=1$, to
\cite{MR1465640,MR2760897,MR2926763} to the Coulomb case in dimension $d=2$
with support restriction on a line, to
\cite
{MR1660943,MR1606719,MR1746976,MR2926763,MR1485778,2012arXiv12013503S,2012arXiv12070718Y}
for the Coulomb case in dimension $d=2$. We also refer to
\cite{2008arXiv08124224B} for the asymptotic analysis in terms of large
deviations of Coulomb determinantal point processes on compact
manifolds of
arbitrary dimension.

The asymptotic analysis of $\mu_N$ as $N\to\infty$ for Coulomb
interactions in
dimension $d\geq3$ motivates our next result, which is stated for the more
general Riesz interactions in dimension $d\geq1$. The \emph{Riesz}
interactions correspond to the choice $W(x,y)=k_{\Delta_\alpha}(x-y)$ where\vspace*{1pt}
$k_{\Delta_\alpha}$, $0<\alpha<d$, $d\geq1$, is the \emph{Riesz
kernel} defined on
$\mathbb{R}^d$, by
%
%
\begin{equation}
\label{eqriesz} k_{\Delta_\alpha}(x):=\frac{1}{{{\vert x \vert
}}^{d-\alpha}}.
\end{equation}
Up to a multiplicative constant, this is the fundamental solution of a
fractional Laplace equation [which is the true Laplace equation
(\ref{eqfs-coul}) when $\alpha=2$], namely
%
%
\begin{equation}
\label{eqfs-riesz} \qquad -c_\alpha\Delta_\alpha k_{\Delta_\alpha
}=\mathcal
{F}^{-1}(1)=\delta_0 \qquad\mbox{with } c_\alpha:=
\frac{\pi^{\alpha-(d/{2})}}{4\pi^2} \frac{\Gamma((d-\alpha
)/{2})}{\Gamma({\alpha}/{2})},
\end{equation}
where the Fourier transform $\mathcal{F}$ and the fractional Laplacian
$\Delta_\alpha
$ are
given by
\[
\mathcal{F}(k_{\Delta_\alpha}) (\xi):=\int_{\mathbb{R}^d}
e^{2i\pi\xi
\cdot x} k_{\Delta_\alpha
}(x) \,dx %
=\frac{1}{c_\alpha4\pi^2{{\vert\xi\vert}}^\alpha}
\]
and
\[
\Delta_\alpha f:= -4\pi^2
\mathcal{F}^{-1} \bigl({{\vert\xi\vert}}^\alpha
\mathcal{F}(f) \bigr).
\]
Note that $\Delta_2=\Delta$ while $\Delta_\alpha$ is a nonlocal
integro-differential
operator when $\alpha\neq2$. When $d\geq3$ and $\alpha=2$ then
Riesz interactions
coincide with Coulomb interactions and the constants match. Beware that our
notation differs slightly from those of Landkof \cite{MR0350027},
page~44. Several
aspects of the Gauss problem in the Riesz case are studied in
\cite{MR2276529,MR2075132,MR2073873}.

In the Riesz case, $0<\alpha<d$, one associates to any probability measure
$\mu$
on~$\mathbb{R}^d$ a function $U_\alpha^\mu\dvtx\mathbb{R}^d\mapsto
[0,+\infty]$
called the
potential of $\mu$ as follows:
\[
U_\alpha^\mu(x):= (k_{\Delta_\alpha}*\mu) (x):=\int
k_{\Delta_\alpha}(x-y) \,d\mu(y)\qquad\forall x\in\R^d.
\]
We refer to Section~\ref{secpotential} for a review of basic
definitions from
potential theory. In particular, one defines there a notion of capacity of
sets, and a property is said to hold quasi-everywhere if it holds
outside a
set of zero capacity. The following theorem is essentially the analogue in
$\mathbb{R}^d$ of a result of Dragnev and Saff on spheres~\cite{MR2276529}. The
analogue problem on compact subsets, without external field, was initially
studied by Frostman \cite{frostman}; see also the book of Landkof
\cite{MR0350027}. A confinement (by an external field or by a support
constraint) is always needed for such type of results.

%
%
\begin{theorem}[(Riesz gases)]\label{thriesz}
Suppose that $W$ is the Riesz kernel\break  $W(x,y)= k_{\Delta_\alpha
}(x-y)$. Then:
\begin{longlist}[(7)]
\item[(1)] The functional $I$ is strictly convex where it is finite.
\item[(2)]\textup{(H1)--(H4)} are satisfied, and Theorem~\ref{thldp}
applies.
\item[(3)] There exists a unique $\mu_\star\in\mathcal{M}_1(\mathbb{R}^d)$
such that
\[
I(\mu_\star)=\inf_{\mu\in\mathcal{M}_1(\mathbb{R}^d)}I(\mu).
\]
\item[(4)] If we define $(\mu_N)_N$ on a unique probability space [for a
sequence $\beta_N\gg N\log(N)$], then with
probability one,
\[
\lim_{N\to\infty}\mu_N=\mu_\star.
\]
If we denote by $C_\star$ the real number
\[
C_\star%
= \int{{ \bigl(U_\alpha^{\mu_\star} + V \bigr)}} \,d
\mu_\star%
= J(\mu_\star) + \int V\,d\mu_\star,
\]
%
then the following additional properties hold:

\item[(5)] The minimizer $\mu_\star$ has compact support, and satisfies
%
%
\begin{eqnarray}
\label{eqrob1} U_\alpha^{\mu_\star}(x) + V(x) &\geq&
C_\star\qquad\mbox{quasi-everywhere},
\\
\label{eqrob2} U_\alpha^{\mu_\star}(x) + V(x) &=& C_\star
\qquad\mbox{for all }x\in\operatorname{supp}(\mu_\star).
\end{eqnarray}
\item[(6)] If a compactly supported measure $\mu$ creates a potential
$U_\alpha^\mu$ such that,
for some constant $C\in\mathbb{R}$,
%
%
\begin{eqnarray}
\label{eqrob3} U_\alpha^\mu(x) + V(x) &=& C\qquad\mbox{on }
\operatorname{supp}(\mu),
\\
\label{eqrob4} U_\alpha^\mu+ V &\geq& C\qquad\mbox{quasi-everywhere},
\end{eqnarray}
then $C = C_\star$ and $\mu=\mu_\star$. The same is true under the
weaker assumptions
%
%
\begin{eqnarray}
\label{eqrob3prime} U_\alpha^\mu(x) + V(x) &\leq& C\qquad\mbox{on }
\operatorname{supp}(\mu),
\\
\label{eqrob4prime} U_\alpha^\mu+ V &\geq& C\qquad\mbox{q.e. on }
\operatorname{supp}(\mu_\star).
\end{eqnarray}
\item[(7)] If $\alpha\leq2$, for any measure $\mu$, the following
``converse'' to
(\ref{eqrob3prime}), (\ref{eqrob4prime}) holds:
%
%
\begin{eqnarray}
\label{eqrob2prime} \sup_{\operatorname{supp}(\mu)} {{ \bigl
(U_\alpha
^\mu+
V \bigr)}} &\geq& C_\star,
\\
\label{eqrob1prime} \mbox{``}\inf_{\operatorname{supp}(\mu
_\star)}\mbox{''}
{{ \bigl(U_\alpha^\mu(x) + V(x) \bigr)}} &\leq&
C_\star,
\end{eqnarray}
where the ``$\inf$'' means that the infimum is taken quasi-everywhere.
\end{longlist}
\end{theorem}

The constant $C_\star$ is called the ``modified Robin constant'' (see,
e.g., \cite{MR1485778}), where the properties
(\ref{eqrob1})--(\ref{eqrob2}) and the characterization
(\ref{eqrob3})--(\ref{eqrob4}) are established for the logarithmic
potential in dimension $2$. The minimizer $\mu_\star$ is called the
\emph{equilibrium measure}.
%
%
\begin{corollary}[(Equilibrium of Coulomb gases with radial external
fields in dimension $\geq3$)]\label{cocoulomb}
Suppose that for a fixed real parameter $\beta>0$, and for every
$x,y\in
\mathbb{R}^d$,
$d\geq3$,
\[
V(x)=v \bigl({{\vert x \vert}} \bigr)\quad\mbox{and}\quad
W(x,y)=\beta
k_{\Delta}(x-y),
\]
where $v$ is two times differentiable. Denote by $d\sigma_r$ the Lebesgue
measure on the sphere of radius $r$, and let $\sigma_d$ be the total
mass of
$d\sigma_1$ (i.e., the surface of the unit sphere of $\mathbb{R}^d$). Let
$w(r) =
r^{d-1}v'(r)$, and suppose either that $v$ is convex, or that $w$ is
increasing. Define two radii $r_0<R_0$ by
\[
r_0 = \inf{{ \bigl\{r>0; v'(r)>0 \bigr\}}} \quad
\mbox{and}\quad w(R_0) = \beta(d-2).
\]
Then the equilibrium measure $\mu_\star$ is supported on the
ring ${{ \{x; {{\vert x \vert}}\in[r_0,R_0] \}}}$ and
is absolutely
continuous with respect to Lebesgue measure
\[
d\mu(r) = M(r) \,d\sigma_r \,dr
\qquad\mbox{where } M(r) =
\frac{w'(r)}{\beta(d-2)\sigma_d r^{d-1}} \IND_{[r_0,R_0]}(r).
\]
In particular, when $v(t)=t^2$, then $\mu_\star$ is the uniform distribution
on the centered ball of radius
\[
{{ \biggl(\beta\frac{d-2}{2} \biggr)}}^{1/d}.
\]
\end{corollary}

The result provided by Corollary~\ref{cocoulomb} on Coulomb gases
with radial
external fields
can be found, for instance, in \cite{MR2647570}, Proposition~2.13. It follows
quickly from the Gauss averaging principle and the characterization
(\ref{eqrob3})--(\ref{eqrob4}). For the sake of completeness, we give
a (short)
proof in Section~\ref{ssradialcoulomb}. By using Theorem~\ref{thriesz} with
$\alpha=2$ together with Corollary~\ref{cocoulomb}, we obtain that
the empirical
measure of a Coulomb gas with quadratic external field in dimension
$d\geq3$
tends almost surely to the uniform distribution on a ball when $N\to
\infty$.
This phenomenon is the analogue in arbitrary dimension $d\geq3$ of the
well-known result in dimension $d=2$ for the logarithmic potential with
quadratic
radial external field (where the uniform law on the disc or ``circular law''
appears as a limit for the complex Ginibre ensemble; see, for instance,
\cite{MR1660943,MR1606719}). The study of the equilibrium measure for Coulomb
interaction with nonradially symmetric external fields was initiated recently
in dimension $d=2$ by Bleher and Kuijlaars in a beautiful work
\cite{MR2921180} by using orthogonal polynomials.

The following proposition shows that in the Riesz case, it is possible to
construct a good confinement potential $V$ so that the equilibrium
measure is
prescribed in advance.

%
%
\begin{corollary}[(Riesz gases: External field for prescribed
equilibrium measure)]\label{coprescription}
Let $0<\alpha<d$, $d\geq1$, and $W(x,y):=k_{\Delta_\alpha}$. Let
$\mu_\star$
be a
probability measure with a compactly supported density $f_\star
\in\mathbf{L}^p(\R^d)$ for some $p>d/\alpha$. Define the
confinement potential
\[
V(x):= -U_\alpha^{\mu_\star}(x) + \bigl[|x|^2-R \bigr]_+,
\qquad x\in\R^d,
\]
where $U_\alpha^{\mu_\star}$ is the Riesz potential created by $\mu
_\star$
and $R>0$ is such that $\operatorname{supp}(\mu_\star)\subset B(0,R)$.
Then the
couple of functions $(V,W)$ satisfy \textup{(H1)--(H4)}, and the
functional
\[
\mu\in\mathcal{M}_1 \bigl(\mathbb{R}^d \bigr)\mapsto%
I(\mu):=\int V \,d\mu%
+ \frac{1}{2}\int\!\!\!\int
k_{\Delta_\alpha}(x-y) \,d\mu(x)\,d\mu(y) \in\mathbb{R}\cup\{
+\infty\}
\]
admits $\mu_\star$ as unique minimizer. In particular, the probability
$\mu_\star$ is the almost sure limit of the sequence ${(\mu_N)}_{N}$
(constructed on the same probability space), as soon as $\beta_N\gg
N\log(N)$.
\end{corollary}

\subsection{\texorpdfstring{Outline of the article.}{Outline of the article}}

In the remainder of this introduction (Section~\ref{sscomments-extensions-related}), we give several comments on our
results, their links with different domains, and possible directions for
further research. Section~\ref{secldp} provides the proof of
Theorem~\ref{thldp} (large deviations principle).
Section~\ref{secminimizing} provides the proof of Theorem~\ref{thriesz},
Corollaries~\ref{cocoulomb}~and~\ref{coprescription}. These
proofs rely on several concepts and tools from Potential Theory, which we
recall synthetically and discuss in Section~\ref{secpotential} for
the sake
of clarity and completeness.

\subsection{\texorpdfstring{Comments, possible extensions and related topics.}{Comments, possible extensions and related topics}}\label{sscomments-extensions-related}
\subsubsection{\texorpdfstring{Noncompactly supported equilibrium measures.}{Noncompactly supported equilibrium measures}}

The assumptions made on the external field $V$ in Theorems~\ref{thldp}~and~\ref{thriesz} explain why the equilibrium measure $\mu_\star
$ is
compactly supported. If one allows a weaker behavior of $V$ at
infinity, then
one may produce equilibrium measures $\mu_\star$ which are not compactly
supported (and may even be heavy tailed). This requires that we adapt
some of the
arguments, and one may use compactification as in \cite{MR2926763}.
This might
allow to extend Corollary~\ref{coprescription} beyond the compactly supported
case.

\subsubsection{\texorpdfstring{Equilibrium measure for Riesz interaction with radial external field.}{Equilibrium measure for Riesz interaction with radial external field}}
To the knowledge of the authors, the computation of the equilibrium measure
for Riesz interactions with radial external field, beyond the more specific
Coulomb case of Corollary~\ref{cocoulomb}, is an open
problem, due to the lack of the Gauss averaging principle when $\alpha
\neq2$.

\subsubsection{\texorpdfstring{Beyond the Riesz and Coulomb interactions.}{Beyond the Riesz and Coulomb interactions}}
Theorem~\ref{thriesz} concerns the minimization of the Riesz interaction
potential with an external field $V$, and includes the Coulomb
interaction if
$d\geq3$. In classical Physics, the problem of minimization of the Coulomb
interaction energy with an external field is known as the Gauss variational
problem \cite{frostman,MR0350027,MR2075132,MR2073873}. Beyond the
Riesz and
Coulomb potentials, the driving structural idea behind Theorem~\ref{thriesz}
is that if $W$ is of the form $W(x,y)=k_D(x-y)$ where $k_D$ is the fundamental
solution of an equation $-Dk_D=\delta_0$ for a local differential
operator $D$
such as $\Delta_\alpha$ with $\alpha=2$, and if $V$ is
super-harmonic for $D$,
that is,
$DV\geq0$, then the density of $\mu_\star$ is roughly given by $DV$
up to
support constraints. This can be easily understood formally with Lagrange
multipliers. The limiting measure $\mu_\star$ depends on $V$ and $W$,
and is
thus nonuniversal in general.

\subsubsection{\texorpdfstring{Second-order asymptotic analysis.}{Second-order asymptotic analysis}}
The asymptotic
analysis of
$\mu_N-\mu_\star$ as $N\to\infty$ is a natural problem, which can
be studied
on various classes of tests functions. It is well known that a repulsive
interaction may affect dramatically the speed of convergence, and make it
dependent over the regularity of the test function. In another
direction, one
may take $\beta_N=\beta N^2$ and study the low temperature regime
$\beta\to
\infty$
at fixed $N$. In the Coulomb case, this leads to Fekete points. We
refer to
\cite{2012arXiv12013503S,MR3046995,sandier-serfaty-1d} for the
analysis of the second order when both $\beta\to\infty$ and $N\to
\infty$.
In the one-dimensional case, another type of local universality inside the
limiting support is available in~\cite{2012arXiv12050671G}.

\subsubsection{\texorpdfstring{Edge behavior.}{Edge behavior}}
Suppose that $V$ is radially symmetric
and that
$\mu_\star$ is supported in the centered ball of radius $r$, like in Corollary
\ref{cocoulomb}. Then one may ask if the radius of the particle system
$\max_{1\leq k\leq n}{{\vert x_k \vert}}$ converges to the edge
$r$ of the limiting
support as $N\to\infty$. This is not provided by the weak convergence of~$\mu_N$. The next question is the fluctuation. In the two-dimensional Coulomb
case, a universality result is available for a class of external fields in
\cite{chapec}.

\subsubsection{\texorpdfstring{Topology.}{Topology}}
It is known that the weak topology can be upgraded
to a Wasserstein topology in the classical Sanov theorem for empirical
measures of i.i.d. sequences (see \cite{MR2593592}), provided that
tails are
strong exponentially integrable. It is then quite natural to ask about
such an
upgrade for Theorem~\ref{thldp}.

\subsubsection{\texorpdfstring{Connection to random matrices.}{Connection to random matrices}}
Our initial inspiration came when writing the survey \cite{MR2908617}, from
the role played by the logarithmic potential in the analysis of the Ginibre
ensemble. When $d=2$, $\beta_N=N^2$, $V(x)={{\vert x \vert}}^2$
and $W(x,y)=\beta
k_\Delta(x-y)=\beta\log\frac{1}{{{\vert x-y \vert}}}$ with
$\beta=2$, then $P_N$ is the
law of
the (complex) eigenvalues of the complex Ginibre ensemble
\[
dP_N(x)=Z_N^{-1}e^{-N\sum_{i=1}^N{{\vert x_i \vert}}^2}\prod
_{i<j}{{\vert x_i-x_j \vert
}}^2\,dx
\]
(here $\mathbb{R}^2\equiv\mathbb{C}$ and $P_N$ is the law of the
eigenvalues of a random
$N\times N$ matrix with i.i.d. complex Gaussian entries of covariance
$\frac{1}{2N}I_2$). For a nonquadratic $V$, we may see $P_N$ as the
law of
the spectrum of random normal matrices such as the ones studied in
\cite{MR2817648}. On the other hand, in the case where $d=1$ and
$V(x)={{\vert x \vert}}^2$ and $W(x,y)=\beta\log\frac
{1}{{{\vert x-y \vert}}}$ with $\beta>0$, then
\[
dP_N(x)=Z_N^{-1}e^{-N\sum_{i=1}^N{{\vert x_i \vert}}^2}\prod
_{i<j}{{\vert x_i-x_j \vert
}}^\beta\,dx.
\]
This is known as the $\beta$-Ensemble in Random Matrix Theory. For
$\beta
=1$, we
recover the law of the eigenvalues of the Gaussian orthogonal ensemble (GOE)
of random symmetric matrices, while for $\beta=2$, we recover the law
of the
eigenvalues of the Gaussian Unitary Ensemble (GUE) of random Hermitian
matrices. It is worthwhile to notice that $-\log{{\vert\cdot
\vert}}$ is the Coulomb
potential in dimension $d=2$, and not in dimension $d=1$. For this
reason, we
may interpret the eigenvalues of GOE/GUE as being a system of charged
particles in dimension $d=2$, experiencing Coulomb repulsion and an external
quadratic field, but constrained to stay on the real axis. We believe
this type of support constraint can be incorporated in our initial
model, at
the price of slightly heavier notation and analysis.

\subsubsection{\texorpdfstring{Simulation problem and numerical approximation of the equilibrium measure.}{Simulation problem and numerical approximation of the equilibrium measure}}
It is natural to ask about the best way to simulate
the probability measure $P_N$. A pure rejection algorithm is too naive. Some
exact algorithms are \mbox{available} in the determinantal case $d=2$ and
$W(x,y)=-2\log{{\vert x-y \vert}}$; see \cite{MR2216966},
Algorithm~18 and
\cite{Scardicchioetal09}. One may prefer to use a nonexact algorithm such
as a Hastings--Metropolis algorithm. One may also use an Euler scheme to
simulate a stochastic process for which $P_N$ is invariant, or use a
Metropolis adjusted Langevin approach (MALA) \cite{MR1888450}. In this
context, a very natural way to approximate numerically the equilibrium measure
$\mu_\star$ is to use a simulated annealing stochastic algorithm.

\subsubsection{\texorpdfstring{More general energies.}{More general energies}}
\hspace*{-1pt}The density of $P_N$ takes the form\break  $\prod_{i=1}^Nf_1(x_i)\* \prod
_{1\leq
i<j\leq
N}f_2(x_i,x_j)$, which\vspace*{1pt} comes from the structure of $H_N$. One may
study more
general energies with many bodies interactions, of the form, for some
prescribed symmetric $W_k\dvtx(\mathbb{R}^d)^k\mapsto\mathbb{R}$,
$1\leq k\leq K$,
$K\geq1$,
\[
H_N(x_1,\ldots,x_N)%
=\sum
_{k=1}^K \sum
_{i_1<\cdots<i_k}N^{-k}W_k(x_{i_1},\ldots,x_{i_k}).
\]
This leads to the following candidate for the asymptotic first-order global
energy functional:
\[
\mu\mapsto\sum_{k=1}^K 2^{-k}
\int\cdots\int W_k(x_1,\ldots,x_k) \,d\mu
(x_1)\cdots d\mu(x_k).
\]

\subsubsection{\texorpdfstring{Stochastic processes.}{Stochastic processes}}\label{sssSDE}
Under general assumptions on $V$ and $W$ (see, e.g., \cite{MR2352327}),
the law $P_N$ is the invariant probability measure of a well-defined (the
absence of explosion comes from the assumptions on $V$ and $W$) reversible
Markov diffusion process ${(X_t)}_{t\in\mathbb{R}_+}$ with state space
\[
{{ \bigl\{x\in\bigl(\mathbb{R}^d \bigr)^N\dvtx H_N(x)<
\infty\bigr\}}} ={{ \biggl\{x\in\bigl(\mathbb{R}^d \bigr)^N\dvtx
\sum_{i<j}W(x_i,x_j)<\infty
\biggr\}}},
\]
solution of the system of Kolmogorov stochastic differential equations
\[
dX_t=\sqrt{2\frac{\alpha_N}{\beta_N}} \,dB_t-
\alpha_N\nabla H_N(X_t) \,dt,
\]
where ${(B_t)}_{t\geq0}$ is a standard Brownian motion on $(\mathbb
{R}^d)^N$ and
where $\alpha_N>0$ is an arbitrary scale parameter (natural choices being
$\alpha_N=1$ and $\alpha_N=\beta_N$). The law $P_N$ is the equilibrium
distribution of
a system of $N$ interacting Brownian particles
${(X_{1,t})}_{t\geq0},\ldots,{(X_{N,t})}_{t\geq0}$ in $\mathbb
{R}^d$ at inverse\vspace*{1pt}
temperature $\beta_N$, with equal individual ``charge'' $q_N:=1/N$,
subject to a
confining potential $\alpha_N V$ acting on each particle and to an interaction
potential $\alpha_N W$ acting on each pair of particles, and one can
rewrite the
stochastic differential equation above as the system of coupled stochastic
differential equations ($1\leq i\leq N$)
\[
dX_{i,t} %
=\sqrt{2\frac{\alpha_N}{\beta_N}} \,dB_{i,t}
-q_N\alpha_N\nabla V(X_{i,t}) -\sum
_{j\neq i}q_N^2
\alpha_N\nabla_1W(X_{i,t},X_{j,t}) \,dt,
\]
where ${(B_t^{(1)})}_{t\geq0},\ldots,{(B_t^{(N)})}_{t\geq0}$ are
i.i.d. standard Brownian motions on $\mathbb{R}^d$. From a partial
differential equations
point of view, the probability measure $P_N$ is the steady state
solution of
the Fokker--Planck evolution equation ${\partial}_t-L=0$ where $L$ is
the elliptic
Markov diffusion operator (second-order linear differential operator without
constant term)
\[
L:=\frac{\alpha_N}{\beta_N}{{ (\Delta-\beta_N\nabla H_N\cdot
\nabla)}},
\]
acting as $Lf=\frac{\alpha_N}{\beta_N}(\Delta f-{{ \langle
\beta_N\nabla H_N,\nabla f \rangle}})$. This
self-adjoint operator in $\rL^2(P_N)$ is the infinitesimal generator
of the
Markov semigroup ${(P_t)}_{t\geq0}$, $P_t(f)(x):=\mathbb{E}(f(X_t)|X_0=x)$.
Let us
take $\alpha_N=\beta_N$ for convenience. In the case where
$V(x)={{\vert x \vert}}^2$ and
$W\equiv0$ (no interaction), then $P_N$ is a standard Gaussian
law $\mathcal{N}(0,I_{dN})$ on $(\mathbb{R}^d)^N$, and
${(X_t)}_{t\geq0}$ is an
Ornstein--Uhlenbeck
Gaussian process; while in the case where $d=1$ and $V(x)={{\vert x
\vert}}^2$ and
$W(x,y)=-\beta\log{{\vert x-y \vert}}$ of some fixed parameter
$\beta>0$, then $P_N$
is the
law of the spectrum of a $\beta$-Ensemble of random matrices, and
${(X_t)}_{t\geq0}$ is a so-called Dyson Brownian motion \cite
{MR2760897}. If
$\mu_{N,t}$ is the law of $X_t$, then $\mathbb{E}\mu_{N,t}\to
\mathbb{E}\mu_N$ weakly as
$t\to\infty$. The study of the dynamic aspects is an interesting problem
connected to McKean--Vlasov models
\cite{MR1440140,MR2062570,li-li-xie,MR3059192,MR3005007}.

\subsubsection{\texorpdfstring{Calogero--(Moser--)Sutherland--Schr\"odinger operators.}{Calogero--(Moser--)Sutherland--Schr\"odinger operators}}
Let us keep the notation used above. We define $U_N:=\beta_NH_N$, and
we take
$\beta_N=N^2$ for simplicity. Let us consider the isometry
$\Theta\dvtx\rL^2(P_N)\to\rL^2(dx)$ defined by
\[
\Theta(f) (x):=f(x)\sqrt{\frac{dP_N(x)}{dx}}=f(x)e^{-({1}/{2})(U_N(x)+\log(Z_N))}.
\]
The differential operator $S:=-\Theta L \Theta^{-1}$ is a Schr\"
odinger operator
\[
S:=-\Theta L \Theta^{-1}=-\Delta+Q, \qquad Q:=\tfrac{1}{4}{{
\vert\nabla U_N \vert}}^2-\tfrac
{1}{2}\Delta
U_N
\]
which acts as $S f=-\Delta f+Qf$. The operator $S$ is self-adjoint in
$\rL^2(dx)$. Being isometrically conjugated, the operators $-L$ and
$S$ have
the same spectrum, and their eigenspaces are isometric. In the case where
$V(x)={{\vert x \vert}}^2$ and $W\equiv0$ (no interactions), we
find that and
$Q=\frac{1}{2}(1-V)$, and $S$ is a harmonic oscillator. On the other hand,
following \cite{MR2641363}, Proposition~11.3.1, in the case $d=1$ and
$W(x,y)=-\log{{\vert x-y \vert}}$ (Coulomb interaction), then $S$
is a
Calogero--(Moser--)Sutherland--Schr\"odinger operator,
\begin{eqnarray*}
S&=&-\Delta-E_0+\frac{1}{4}\sum_{i=1}^Nx_i^2
-\frac{1}{2}\sum_{1\leq i<j\leq N}
\frac{1}{(x_i-x_j)^2},
\\
E_0&:=&\frac{N}{2}+\frac{N(N-1)}{2}.
\end{eqnarray*}
More examples are given in \cite{MR2641363}, Proposition~11.3.2,
related to
classical ensembles of random matrices. The study of the spectrum and
eigenfunctions of such operators is a wide subject, connected to Dunkl
operators. These models attracted some attention due to the fact that for
several natural choices of the potentials $V,W$, they are exactly
solvable (or
integrable). We refer to \cite{MR2641363}, Section~11.3.1, \cite{MR1827871}, Section~9.6, \cite{dunkl-book}, Section~2.7 and
references therein.

\section{\texorpdfstring{Proof of the large deviations principle---Theorem~\protect\ref{thldp}.}{Proof of the large deviations principle---Theorem~1.1}}\label{secldp}

The proof of Theorem~\ref{thldp} is split is several steps.

\subsection{\texorpdfstring{A standard reduction.}{A standard reduction}}

To prove Theorem~\ref{thldp}, we will use the following standard
reduction; see, for instance, \cite{MR2571413}, Chapter~4.

%
%
\begin{proposition}[(Standard reduction)]\label{propreductionldp}
Let $(Q_N)_N$ be a sequence of probability measures on some Polish space
$(\mathcal{X},d)$, $(Z_N)_N$ and $(\varepsilon_N)_N$ two sequences of
positive numbers
with $\varepsilon_N \to0$ and $\mathcal{I}\dvtx\mathcal{X}\to\R
\cup\{
+\infty\}$ be a function
bounded from below.
\begin{longlist}[(2)]
\item[(1)] Suppose that the sequence $(Q_N)_N$ satisfies the following conditions:
\begin{enumerate}[(a)]
\item[(a)] The sequence $(Z_NQ_N)_N$ is exponentially tight: for all $L\geq0$
there exists a compact set $K_L \subset\mathcal{X}$ such that
\[
\limsup_{N\to\infty} \varepsilon_N \log
Z_NQ_N{{ (\mathcal{X}\setminus K_L )}}
\leq- L.
\]
\item[(b)] For all $x\in\mathcal{X}$,
\[
\lim_{r\to0} \limsup_{N\to\infty}
\varepsilon_N \log Z_NQ_N \bigl(B(x,r) \bigr)
\leq- \mathcal{I}(x),
\]
where $B(x,r):=\{y\in\mathcal{X}\dvtx d(x,y)\leq r\}$.
\end{enumerate}
Then the sequence $(Z_NQ_N)_N$ satisfies the following large deviation
upper bound: for all Borel set $A \subset\mathcal{X}$, it holds
%
%
\begin{equation}
\label{eqbornesupprov} \limsup_{N\to\infty} \varepsilon_N\log
Z_NQ_N(A) %
\leq-\inf{{ \bigl\{\mathcal{I}(\mu); \mu\in\operatorname
{clo}(A) \bigr
\}}}.
\end{equation}
\item[(2)] If, in addition, $(Z_NQ_N)_N$ satisfies the following large deviation
lower bound: for any Borel set $A\subset\mathcal{X}$,
%
%
\begin{equation}
\label{eqborneinfprov} -\inf\bigl\{\mathcal{I}(x); x\in
\operatorname
{int}(A) \bigr\}
\leq\liminf_{N\to\infty} \varepsilon_N \log
Z_NQ_N(A),
\end{equation}
then $(Q_N)_N$ satisfies the full large deviation principle with speed
$\varepsilon_N$ and rate function $\mathcal{I}_\star= \mathcal
{I}-\inf_{x\in\mathcal{X}} \mathcal{I}(x)$,
namely for any Borel set $A\subset\mathcal{X}$,
\begin{eqnarray*}
-\inf\bigl\{\mathcal{I}_\star(x); x\in\operatorname{int}(A) \bigr
\} &\leq&\liminf_{N\to\infty} \varepsilon_N \log
Q_N(A)
\\
&\leq&\limsup_{N\to\infty} \varepsilon_N\log
Q_N(A)
\\
&\leq& -\inf{{ \bigl\{\mathcal{I}_\star(x); x \in
\operatorname{clo}(A) \bigr\}}}.
\end{eqnarray*}
\end{longlist}
\end{proposition}

\begin{pf}
Let us begin by (1). Let $\delta>0$; by assumption, for any $x\in
\mathcal{X}$,
there is $\eta_x>0$ such that
\[
\limsup_{N\to\infty} \varepsilon_N \log
Z_NQ_N \bigl(B(x,\eta_x) \bigr) \leq-
\mathcal{I}(x) + \delta.
\]
If $F\subset\mathcal{X}$ is compact, there is a finite family
$(x_i)_{1\leq
i\leq
m}$ of points of $F$ such that $F \subset\bigcup_{i=1}^m B(x_i,\eta_{x_i})$.
Therefore,
\begin{eqnarray*}
\limsup_{N\to\infty}\varepsilon_N\log
Z_NQ_N(F) &\leq&\limsup_{N\to\infty
}
\varepsilon_N\log\Biggl( \sum_{i=1}^NZ_NQ_N
\bigl(B(x_i,\eta_{x_i}) \bigr) \Biggr)
\\
&=& \max_{1\leq i\leq m} \limsup_{N\to\infty}
\varepsilon_N\log\bigl( Z_NQ_N
\bigl(B(x_i,\eta_{x_i}) \bigr) \bigr)
\\
&\leq&\max_{1\leq i\leq m} -\mathcal{I}(x_i)+\delta
\\
&\leq& -\inf_F \mathcal{I}+ \delta.
\end{eqnarray*}
Letting $\delta\to0$ yields to (\ref{eqbornesupprov}) for $A=F$ compact.

Now if $F$ is an arbitrary closed set, then for all $L>0$, since $F\cap
\mathcal{K}_L$ is compact, it holds
\begin{eqnarray*}
&&\limsup_{N\to\infty} \varepsilon_N\log
Z_NQ_N(F)
\\
&&\qquad\leq\max\Bigl( \limsup_{N\to\infty} \varepsilon_N
\log Z_N Q_N(F\cap K_L), \limsup
_{N\to\infty} \varepsilon_N\log Z_NQ_N
\bigl(K_L^c \bigr) \Bigr)
\\
&&\qquad\leq\max\Bigl(-\inf_{F\cap K_L} \mathcal{I}; -L \Bigr).
\end{eqnarray*}
Letting $L\to\infty$ shows that (\ref{eqbornesupprov}) is true for
arbitrary closed sets $F$. Since $A \subset\operatorname{clo}(A)$,
the upper bound
(\ref{eqbornesupprov}) holds for arbitrary Borel sets $A$.

To prove (2), take $A=\mathcal{X}$ in (\ref{eqborneinfprov}) and
(\ref{eqbornesupprov}) to get
\[
\lim_{N\to\infty} \varepsilon_N \log(Z_N) =
-\inf\mathcal{I}\in\R.
\]
Subtracting this to (\ref{eqborneinfprov}) and (\ref{eqbornesupprov})
gives the large deviations principle with rate function $\mathcal
{I}_\star$.
\end{pf}

In our context, $\mathcal{X}=\mathcal{M}_1(\mathbb{R}^d)$ is
equipped with
the Fortet--Mourier
distance~(\ref{eqFortetMourier}).

\subsection{\texorpdfstring{Properties of the rate function.}{Properties of the rate function}}
In the following lemma, we prove different properties of the rate function
$I_\star$ including those announced in Theorem~\ref{thldp}, point~(1).

%
%
\begin{lemma}[(Properties of the rate function)]\label{lemgoodrate}
Under assumptions \textup{(H1)--(H3)}:
\begin{longlist}[(4)]
\item[(1)] $I_\star$ is well defined;
\item[(2)] $I_\star(\mu)<\infty$ implies $\int|V| \,d\mu<\infty$ and
$\int\!\!\int
|W| \,d\mu^2<\infty$;
\item[(3)] $I_\star(\mu)<\infty$ for any compactly supported probability
$\mu$
with a bounded density with respect to Lebesgue;
\item[(4)] $I_\star$ has is a good rate function (i.e., the levels sets
$\{I_\star\leq k\}$ are compact).
\end{longlist}
\end{lemma}

\begin{pf}
Let us define $\varphi\dvtx\mathbb{R}^d\times\mathbb{R}^d\to
(-\infty,+\infty
]$ by
$\varphi(x,y):=\frac{1}{2} (V(x)+V(y)+W(x,y) )$.
\begin{longlist}[(4)]
\item[(1)] Since $V$ is continuous and $V(x)\to\infty$ as $|x|\to\infty$
thanks to \textup{(H2)}, the function $V$ is bounded from below. Using
\textup{(H3)} it follows that $\varphi$ is bounded from below. The functional
$I_\star$ is thus well defined with values in $[0,\infty]$.
\item[(2)] Assume that $I(\mu) = \int\!\!\int\varphi \,d\mu^2 <\infty$.
Since $V$ is
bounded from below,
$[V]_-\in\mathbf{L}^1(\mu)$.
From \textup{(H3)} and the definition of $\varphi$,
\[
2\varphi(x,y) = V(x) + V(y) + W(x,y) \geq c + (1-\varepsilon_0)
\bigl(V(x) + V(y) \bigr).
\]
Therefore
\[
(1-\varepsilon_0) \bigl([V]_+(x) + [V]_+(y) \bigr) \leq2
\varphi(x,y) - c + (1-\varepsilon_0) \bigl([V]_-(x) + [V]_-(y)
\bigr),
\]
so $[V]_+ \in\mathbf{L}^1(\mu)$ and $\int{{\vert V \vert}}
\,d\mu< \infty$. Since
%
%
\begin{equation}
c- \varepsilon_0V(x) - \varepsilon_0 V(y) \leq W(x,y)
\leq2\varphi(x,y) - V(x) - V(y), \label{eqdoubleBoundW}
\end{equation}
this implies that $W\in\mathbf{L}^1(\mu^2)$.
\item[(3)] It is clearly enough to prove that $W$ is locally Lebesgue integrable
on \mbox{$\R^d \times\R^d$}. Let $K$ be a compact of $\R^d$; according to
\textup{(H2)} and \textup{(H3)} the function $W$ is bounded from
below on
$K\times K$. On the other hand, letting
\[
\alpha_K(z)=\sup\bigl\{W(x,y); |z-y|\geq|z|, x,y \in K \bigr\},
\]
we have $W(x,y)\leq\alpha_K(x-y)$, for all $x,y \in K$. Assumption
\textup{(H1)} then easily implies that $(x,y) \mapsto\alpha_K(x-y)$ is
integrable on $K\times K$.
\item[(4)] According to the monotone convergence theorem,
\[
I = \sup_{n\in\mathbb{N}} I_n, \qquad I_n(
\mu):=\int\!\!\!\int\min\bigl( \varphi(x,y); n \bigr) \,d\mu
(x)\,d\mu(y).
\]
The functions $\min(\varphi,n)$ being bounded and continuous, it follows
that the functionals $I_n$ are continuous for the weak topology; see, for
instance, \cite{MR2571413}, Lemma~7.3.12. Being a supremum of continuous
functions, $I$ is lower semi-continuous. Set $b_\star=\inf\varphi$; we
have, for every $\mu\in\mathcal{M}_1(\mathbb{R}^d)$, $L>0$,
\begin{eqnarray*}
I(\mu)-b_\star%
&=& \int\!\!\!\int\bigl(\varphi(x,y)-b_\star
\bigr) \,d\mu(x)\,d\mu(y)
\\
&\geq&\int\!\!\!\int\IND_{{{\vert x \vert}}>L,{{\vert y
\vert}}>L} \bigl(\varphi(x,y)-b_\star
\bigr) \,d\mu(x)\,d\mu(y)
\\
&\geq&(b_L-b_\star)\mu\bigl({{\vert x \vert}}>L
\bigr)^2,
\end{eqnarray*}
where $b_L:= \inf_{|x|>L, |y|>L} \varphi(x,y)$. According to \textup{(H2)}
and \textup{(H3)}, we see that $b_L\to+\infty$ as $L\to+\infty$.
Therefore, there exists $L_\star>0$ such that $b_L>b_\star$ for every
$L>L_\star$. We get then for every real number $r\geq b_\star$,
\[
\bigl\{\mu\in\mathcal{M}_1 \bigl(\mathbb{R}^d \bigr)\dvtx I(\mu
)\leq
r \bigr\} \subset{{ \biggl\{\mu\in\mathcal{M}_1 \bigl(
\mathbb{R}^d \bigr)\dvtx\mu\bigl({{\vert x \vert}}>L \bigr)\leq
\sqrt{
\frac{r-b_\star}{b_L-b_\star}},L>L_\star\biggr\}}}.
\]
Since $b_L\to+\infty$ as $L\to+\infty$, the subset of $\mathcal
{M}_1(\mathbb{R}^d)$ on
the right-hand side is tight, and the Prohorov theorem implies then that
it is relatively compact for the topology of $\mathcal{M}_1(\mathbb{R}^d)$.
Since $I$ is
lower semi-continuous, the set $\{I\leq r\}$ is also closed, which
completes the proof.\quad\qed
\end{longlist}\noqed
\end{pf}

\subsection{\texorpdfstring{Proof of the upper bound.}{Proof of the upper bound}}

For\vspace*{2pt} all $N\geq1$, one denotes by $Q_N$ the law of
$\mu_N=\frac{1}{N}\sum_{i=1}^N \delta_{x_i}$ under the probability $P_N$
defined by (\ref{eqPN}): $Q_N$ is an\vspace*{1pt} element of $\mathcal
{M}_1(\mathcal{M}_1(\mathbb{R}^d))$.

%
%
\begin{lemma}[(Exponential tightness)]\label{lemexptight}
If $\beta_N \gg N$, then under assumptions
\textup{(H2)--(H3)}, the sequence of\vspace*{1pt} measures $(Z_NQ_N)_{N}$ is
exponentially tight: for all $L\geq0$ there exists a compact set $K_L
\subset\mathcal{M}_1(\mathbb{R}^d)$ such that
%
%
\begin{equation}
\label{eqexpTight} \limsup_{N\to\infty} %
\frac{\log Z_NQ_N{{ (\mathcal{M}_1(\mathbb{R}^d)\setminus K_L
)}}}{\beta_N}
\leq- L.
\end{equation}
\end{lemma}

\begin{pf}
For any $L\geq0$, let $L':= \frac{L-c/2}{1-\varepsilon_o}$ and set
$K_L:=\{\mu\in\mathcal{M}_1(\mathbb{R}^d);\break  \int V \,d\mu\leq L'\}$.
Since~\textup{(H2)}
holds, $V(x)\to\infty$ when $|x|\to+\infty$ and $V$ is continuous. By
Prohorov's theorem on tightness, this implies that $K_L$ is compact in
$\mathcal{M}_1(\mathbb{R}^d)$.

It remains to check (\ref{eqexpTight}). Let us consider the law $\nu
_V\in
\mathcal{M}_1(\mathbb{R}^d)$ defined by
%
%
\begin{equation}
\label{eqdefNuV} d\nu_V(x):=\frac{e^{-V(x)}}{C_V} \,dx, %
\qquad
C_V:=\int e^{-V(x)} \,dx>0.
\end{equation}
Using (\ref{eqdoubleBoundW}) to bound $W$ from below, we get
\begin{eqnarray*}
\hspace*{-3pt}&& Z_NQ_N{{ \biggl(\int V \,d\mu_N
>L' \biggr)}}
\\
\hspace*{-3pt}&&\quad=\int_{(\mathbb{R}^d)^N} \IND_{{ \{\int V \,d\mu_N >
L' \}}} \exp{{ \biggl(-
\frac{\beta_N}{2}\int\!\!\!\int_{\neq}W \,d\mu
_N^2-\beta_N\int V \,d\mu_N
\biggr)}} \,dx
\\
\hspace*{-3pt}&&\quad\leq\int_{(\mathbb{R}^d)^N} \IND_{{ \{\int V \,d\mu_N >
L' \}}}
\\
\hspace*{-3pt}&&\hspace*{48pt}{}\times \exp \biggl(-
\frac{\beta_N}{2}\int\!\!\!\int_{\neq
} \bigl(c-\varepsilon
_o \bigl(V(x)+V(y) \bigr) \bigr) \,d\mu_N^2-
\beta_N\int V \,d\mu_N \biggr) \,dx
\\
\hspace*{-3pt}&&\quad=\int_{(\mathbb{R}^d)^N} \IND_{{ \{\int V \,d\mu_N >
L' \}}}
\\
\hspace*{-3pt}&&\hspace*{48pt}{}\times \exp \biggl(-
\frac{\beta_N}{2}c\frac{N-1}{N}-\beta_N{{ \biggl(1-
\varepsilon_o\frac{N-1}{N} \biggr)}}\int V \,d\mu_N
\biggr) \,dx
\\
\hspace*{-3pt}&&\quad=C_V^N\int_{(\mathbb{R}^d)^N}\IND_{{ \{\int V \,d\mu_N >L' \}}}
\\
\hspace*{-3pt}&&\hspace*{66pt}{}\times \exp \biggl(-\frac{\beta_N}{2}c\frac{N-1}{N}
\\
\hspace*{-3pt}&&\hspace*{100pt}{} -
\biggl(\beta_N{{ \biggl(1-\varepsilon_o
\frac{N-1}{N} \biggr)}}-N \biggr)\int V \,d\mu_N \biggr) \,d
\nu_V^{\otimes N}(x).
\end{eqnarray*}
Now, if $N$ is large enough, then $\beta_N{{ (1-\varepsilon
_o\frac{N-1}{N} )}}\geq N$,
so that
\begin{eqnarray*}
&& Z_NQ_N{{ \biggl(\int V \,d\mu_N
>L' \biggr)}}
\\
&&\qquad \leq C_V^N \exp{{ \biggl(-
\frac{\beta_N}{2}c\frac{N-1}{N} \biggr)}} \exp{{ \biggl(-{{
\biggl(
\beta_N{{ \biggl(1-\varepsilon_o\frac
{N-1}{N}
\biggr)}}-N \biggr)}}L' \biggr)}}.
\end{eqnarray*}
Therefore, when $N$ is large enough, using the fact that $\beta_N \gg N$,
\begin{eqnarray*}
&& \frac{\log Z_NQ_N{{ (\int V \,d\mu_N>L' )}}}{\beta_N}
\\
&&\qquad \leq \frac
{N\log
C_V}{\beta_N} -\frac{1}{2}c
\frac{N-1}{N} -{{ \biggl({{ \biggl(1-\varepsilon_o
\frac{N-1}{N} \biggr)}}-\frac
{N}{\beta_N} \biggr)}}L'
\\
&&\qquad =-\frac{1}{2}c-(1-\varepsilon_o)L'+o_{N\to\infty}(1)
\\
&&\qquad =-L + o_{N\to\infty}(1).
\end{eqnarray*}
This implies (\ref{eqexpTight}) and completes the proof.
\end{pf}

%
%
\begin{proposition}[(Upper bound)]\label{propbornesup}
If $\beta_N\gg N$, then under assumptions
\mbox{\textup{(H2)--(H3)}}, for all $r\geq0$, for all $\mu\in\mathcal
{M}_1(\mathbb{R}^d)$,
\[
\lim_{r \to0}\limsup_{N \to+\infty} %
\frac{\log Z_N Q_N(B(\mu,r))}{\beta_N} \leq-I(\mu),
\]
where the ball $B(\mu,r)$ is defined for the Fortet--Mourier distance
(\ref{eqFortetMourier}).
\end{proposition}

\begin{pf}
In contrast with the proof of Lemma~\ref{lemexptight}, our
objective now
is to keep enough empirical terms inside the exponential in order to get
$I(\mu)$ at the limit. Introduce $\varphi(x,y)=\frac
{1}{2}(W(x,y)+V(x)+V(y))$, $x,y\in\R^d$. According to \textup{(H3)},
it holds
%
%
\begin{equation}
\label{eqminorationvphi} \varphi(x,y)\geq\frac{c}{2} + \frac
{1-\varepsilon
_o}{2}
\bigl(V(x)+V(y) \bigr)\qquad\forall x,y\in\R^d,
\end{equation}
for\vspace*{2pt} some $c\in\R$ and $\varepsilon_o\in(0,1)$.
Define $\lambda_N= \frac{N^2}{(1-\varepsilon_o)(N-1)}$, and let us
bound the
function $H_N$ from below using (\ref{eqminorationvphi}) at the third
line: for all $n\in\mathbb{N}$, it holds
\begin{eqnarray*}
\beta_NH_N(x)
&=& \beta_N \biggl(
\frac{1}{2}\int\!\!\!\int_{\neq} W \,d\mu
_N^2 + \int V \,d\mu_N \biggr)
\\
&=&
\beta_N \biggl(\int\!\!\!\int_{\neq} \varphi \,d\mu
_N^2 + \frac
{1}{N}\int V \,d\mu_N
\biggr)
\\
&\geq &(\beta_N-\lambda_N)\int\!\!\!\int
_{\neq} \varphi \,d\mu_N^2 + \lambda
_N \int\!\!\!\int_{\neq} \varphi \,d
\mu_N^2 + \frac{\beta_N}{N}\min V
\\
&\geq& (\beta_N-\lambda_N)\int\!\!\!\int
_{\neq} \varphi \,d\mu_N^2 + \lambda
_N \frac
{(N-1)c}{2N}+N\int V \,d\mu_N + \frac{\beta_N}{N}
\min V
\\
&\geq& (\beta_N-\lambda_N)\int\!\!\!\int\varphi\wedge n \,d\mu_N^2 -(\beta_N-\lambda_N)
\frac{n}{N}
\\
&&{} + \lambda_N \frac{(N-1)c}{2N}+N\int V \,d
\mu_N + \frac
{\beta_N}{N}\min V
\\
&=& (\beta_N-\lambda_N)\int\!\!\!\int\varphi\wedge n \,d
\mu_N^2 +N\int V \,d\mu_N +o(
\beta_N),
\end{eqnarray*}
since $\beta_N\gg N$ and $\lambda_N = O(N)$.

Denoting by $I_n(\nu)=\int\!\!\int\varphi\wedge n \,d\nu^2$, $\nu
\in
\mathcal{M}_1(\mathbb{R}
^d)$, and using the preceding lower bound, we see that for every $\mu
\in\mathcal{M}_1(\mathbb{R}^d)$, $r\geq0$ and
$N\gg1$, we have
\begin{eqnarray*}
&& Z_N Q_N \bigl(B(\mu,r) \bigr)
\\
&&\qquad = \int
_{(\mathbb{R}^d)^N} \IND_{B(\mu,r)}(\mu_N) \exp{{ \bigl(-
\beta_N H_N(x) \bigr)}} \,dx
\\
&&\qquad \leq  e^{o(\beta_N)}\int_{(\mathbb{R}^d)^N} \IND_{B(\mu,r)}(
\mu_N) \exp{{ \bigl(-(\beta_N-\lambda_N)I_n(
\mu_N) \bigr)}} \prod_{i=1}^Ne^{-V(x_i)}\,dx
\\
&&\qquad = C_V^N e^{o(\beta_N)}\int_{(\mathbb{R}^d)^N}
\IND_{B(\mu,r)}(\mu_N) \exp{{ \bigl(- (\beta_N-
\lambda_N)I_n(\mu_N) \bigr)}} \,d
\nu_V^N
\\
&&\qquad \leq C_V^N e^{o(\beta_N)}e^{-(\beta_N-\lambda_N)\inf_{\nu\in
B(\mu,r)}
I_n(\nu)},
\end{eqnarray*}
where the definition of $\nu_V$ is given by (\ref{eqdefNuV}).

Therefore, since $\beta_N\gg N$ and $\lambda_N = O(N)$,
\[
\limsup_{N \to+\infty} \frac{\log Z_N Q_N(B(\mu,r))}{\beta_N} %
\leq- \inf
_{\nu\in B(\mu,r)} I_n(\nu).
\]
Since $\varphi\wedge n$ is bounded continuous, the functional $I_n$ is
continuous for the weak topology. As a result, it holds
\[
\lim_{r\to0} \inf_{\nu\in B(\mu,r)} I_n(
\nu) = I_n(\mu).
\]
Finally, the monotone convergence theorem implies that
$\sup_{n\geq1}I_n(\mu)=I(\mu)$, which ends the proof.
\end{pf}

Using this proposition, Lemma~\ref{lemexptight} and the
first point of Proposition~\ref{propreductionldp}, we get the
upper bound of Theorem~\ref{thldp}, point~(2).

\subsection{\texorpdfstring{The lower bound and the full LDP.}{The lower bound and the full LDP}}
In what follows, we denote by $|A|$ the Lebesgue measure of a Borel set
$A \subset\R^n$.

%
%
\begin{proposition}[(Lower bound for regular probabilities)]\label
{propborneinf}
Under the assumptions \textup{(H1)--(H3)}, if $\beta_N \gg N\log
(N)$, then
for every probability measure $\mu$ on $\R^d$ supported in a box
$B=\prod_{i=1}^d[a_i,b_i]$, $a_i,b_i\in\R$, with a density $h$ with respect
to the Lebesgue measure such that, for some $\delta>0$, $\delta\leq
h\leq
\delta^{-1}$ on $B$, it holds
\[
\liminf_{N\to\infty}\frac{\log Z_NQ_N(B(\mu,r))}{\beta_N}\geq
-I(\mu
)\qquad\forall r
\geq0,
\]
where $B(\mu,r)$ is the open ball of radius $r$ centered at $\mu$ for the
Fortet--Mourier distance~(\ref{eqFortetMourier}).
\end{proposition}

If $B$ is the box $\prod_{k=1}^d [a_k,b_k]$ in $\mathbb{R}^d$, let
$l(B)$ and $L(B)$ be the minimum (resp., maximum) edge length
\[
l(B) = \min_{1\leq k\leq d}( b_k - a_k),
\qquad%
L(B) = \max_{1\leq k\leq d}(
b_k - a_k).
\]
We admit for a moment the following result:

%
%
\begin{lemma}[(Existence of nice partitions)]\label{lemnicepartitions}
For all $d$ and all $\delta>0$ there exists a constant $C(d,\delta)$ such
that the following holds. For any\vspace*{1pt} box $B$, any integer $n$, and any measure
$\mu$ with a density $h$ w.r.t. Lebesgue measure, if $\delta\leq h
\leq
\delta^{-1}$, then there exists a partition $(B_1, B_2, \ldots, B_n)$
of $B$
in $n$ sub-boxes, such that:
\begin{longlist}[(2)]
\item[(1)] $B$ is split in equal parts: for all $i$, $\mu(B_i) =
\frac{1}{n}\mu(B)$;
\item[(2)] the edge lengths of the $B_i$ are controlled
\[
\frac{1}{C(d,\delta) n^{1/d}}l(B) %
\leq l(B_i) %
\leq
L(B_i) %
\leq\frac{C(d,\delta)}{n^{1/d}} L(B).
\]
\end{longlist}
\end{lemma}

\begin{pf*}{Proof of Proposition~\ref{propborneinf}}
For each $N$ we apply Lemma~\ref{lemnicepartitions} to obtain a partition
of $B$ in $N$ boxes $B^N_1, \ldots, B^N_N$. Let $d_N$ be the maximum
diameter of the boxes: by the lemma, since $\mu(B) = 1$,
\[
\frac{c_1}{N^{1/d}} \leq l \bigl(B^N_i \bigr) \quad
\mbox{and}\quad d_N:=\max_{1\leq i\leq N}\sup
_{x,y\in B_i^N}|x-y| \leq\frac{c_2}{N^{1/d}},
\]
where $c_1$ and $c_2$ only depend on $B$, $d$ and $\delta$.

Note that, for all $1$-Lipschitz function $f$ with ${{\llVert f\rrVert
}}_\infty\leq
1$, if
$x_i \in B_i^N$ for all $i\leq N$, since $\mu(B_i) = 1/N$ we have
\begin{eqnarray*}
{{ \Biggl\vert\frac{1}{N}\sum_{i=1}^Nf(x_i)-
\int f \,d\mu\Biggr\vert}} &\leq&\sum_{i=1}^N
\int_{B_i^N} \bigl|f(x)-f(x_i)\bigr| \,d\mu(x)
\\
&\leq& d_N.
\end{eqnarray*}
If $N$ is large enough, $d_N \leq r$, which implies that
\[
{{ \bigl\{(x_1, \ldots, x_n) \in B_1^N
\times\cdots\times B_N^N \bigr\} }} \subset\bigl\{\mu
_N \in B(\mu,r) \bigr\}.
\]
Let us denote by $C_i^N\subset B_i^N$ the box obtained from $B_i^N$ by an
homothetic transformation of center the center of $B_i^N$ and ratio (say)
$1/2$. It holds
\[
Z_NQ_N \bigl(B(\mu,r) \bigr) %
\geq\exp{{
\Biggl(- \frac{\beta_N}{N}\sum_{i=1}^N
\max_{C_i^N}V -\frac{\beta_N}{N^2}\sum
_{i<j}\max_{C_i^N\times C_j^N}W \Biggr)}} %
\prod_{i=1}^N \bigl|C_i^N\bigr|.
\]
Since $ {{\vert C_i^N \vert}}\geq(l(B_i^N)/2)^d \geq c_3/N$ for some
absolute constant $c_3$, we have
\[
\frac{\log\prod_{i=1}^N|C_i^N|}{\beta_N} \geq\frac{N\log
(c_3)}{\beta
_N}-\frac{N\log(N)}{\beta_N} \mathop{\longrightarrow}_{N\to\infty
} 0,
\]
and thus we conclude that
\begin{eqnarray*}
&& \liminf_{N\to+\infty} \frac{\log{{ (Z_NQ_N(B(\mu,r))
)}}}{\beta_N}
\\
&&\qquad \geq- \limsup
_{N\to\infty}\frac{1}{N}\sum_{i=1}^N
\max_{C_i^N}V - \limsup_{N\to\infty}
\frac{1}{N^2} \sum_{i<j}\max_{C_i^N\times C_j^N}W.
\end{eqnarray*}
For all $N$, consider the locally constants functions $V_N\dvtx B\to\R
$ and
$W_N\dvtx B\times B \to\R$ defined by
\[
\forall x\in B_i^N\qquad V_N(x):= \max
_{C_i^N} V
\]
and
\[
\forall(x,y)\in B_i^N\times B_j^N\qquad
W_N(x,y):=\max_{C_i^N\times C_j^N}W.
\]
Since $\mu(B^N_i) = 1/N$, it holds
\[
\frac{1}{N}\sum_{i=1}^N \max
_{C_i^N}V = \int_{B} V_N(x) \,d
\mu(x)
\]
and
\[
\frac{1}{N^2}\sum_{i<j}
\max_{C_i^N\times C_j^N}W = \frac{1}{2}\int_{x\neq y}
W_N(x,y) \,d\mu(x)\,d\mu(y).
\]
The uniform continuity of $V$ on $B$ immediately implies that $V_N$
converges uniformly to $V$, and so
\[
\int V_N \,d\mu\to\int V \,d\mu.
\]
For the same reason $W_N$ converges uniformly to $W$ on
\[
(B\times B)\cap\bigl\{(x,y)\in\R^d\times\R^d;{{\vert
x-y \vert}}\geq u \bigr\},
\]
for all $u>0$. According to \textup{(H2)} and \textup{(H3)}, the function
$W$ is bounded from below on $B\times B$. It follows that the functions
$W_N$ are bounded from below by some constant independent on $N$. To apply
the dominated convergence theorem, it remains to bound $W_N$ from
above by
some integrable function. Let
\[
\alpha_B(u):=\sup_{{{\vert x-y \vert}}\geq u} W(x,y),
\]
so that $W(x,y)\leq\alpha_B({{\vert x-y \vert}})$. Obviously
\[
\max_{(x,y)\in B_i^N\times B_j^N}{{\vert x-y \vert}} \leq2d_N +
\min_{(x,y)\in C_i^N\times C_j^N}{{\vert x-y \vert}}.
\]
By construction, since $i\neq j$, we have
\[
\min_{(x,y)\in C_i^N\times C_j^N}{{\vert x-y \vert}} %
\geq
\frac{1}{4} \bigl(l \bigl(B_i^N \bigr)+l
\bigl(B_j^N \bigr) \bigr) %
\geq
\frac{c_1}{4}N^{-1/d} %
\geq\frac{c_1}{4c_2}d_N.
\]
Therefore,
there is an absolute constant $c_4$ such that
\[
\min_{(x,y)\in C_i^N\times C_j^N}{{\vert x-y \vert}} %
\geq
c_4 \max_{(x,y)\in B_i^N\times B_j^N}{{\vert x-y \vert}}.
\]
Since the function $\alpha_B$ is nonincreasing, it holds
\[
\max_{C_i^N\times C_j^N}\alpha_B \bigl({{\vert x-y \vert
}} \bigr) \leq\min_{(x,y)\in B_i^N\times B_j^N}\alpha_B
\bigl(c_4{{\vert x-y \vert}} \bigr).
\]
We conclude from this that $W_N(x,y)\leq\alpha_B(c_4{{\vert x-y
\vert}})$,
$x\neq y$.
It follows from assumption \textup{(H1)} that the function
$\alpha_B(c_4|x-y|)$ is integrable on $B\times B$ with respect to Lebesgue
measure. Since the density of $\mu$ with respect\vspace*{1pt} to Lebesgue is
bounded from
above this function is integrable on $B\times B$ with respect to $\mu^2$.
Applying the dominated convergence theorem, we conclude that
\begin{eqnarray*}
\liminf_{N\to\infty} \frac{\log{{ (Z_NQ_N(B(\mu,r))
)}}}{\beta_N} &\geq&-\int V(x) \,d\mu(x) -
\frac{1}{2} \int\!\!\!\int W(x,y) \,d\mu(x)\,d\mu(y)
\\
&=& -I(\mu).
\end{eqnarray*}\upqed
\end{pf*}

Let us now prove that ``nice'' partitions exist.

\begin{pf*}{Proof of Lemma~\ref{lemnicepartitions}}
The proof is an induction on the dimension $d$.

\textit{Base case}. Let $d=1$, and suppose that $B = [a_0,b_0]$. Since
$\mu$
has a density, there exist ``quantiles'' $a_0 = q_0<q_1 <\cdots<q_n = b_0$
such that
\[
\forall1\leq i \leq n \qquad\mu\bigl([q_{i-1},q_i]
\bigr) = \frac{1}{n} \mu(B).
\]
In this simple case $l(B) = L(B) = b_0 - a_0$ and $l(B_i) = L(B_i) =
(q_i - q_{i-1})$.
The boundedness assumption on $h$ implies that
\begin{eqnarray*}
\delta(q_i - q_{i-1}) &\leq&\mu\bigl([q_{i-1},q_i]
\bigr) \leq\frac{1}{\delta} (q_i - q_{i-1}),
\\
\delta(b_0 - a_0) &\leq&\mu(B) \leq\frac{1}{\delta}
(b_0 - a_0),
\end{eqnarray*}
and the claim holds for $d=1$ with $C(1,\delta) = 1/\delta^2$.

\textit{Induction step}. Suppose that the statement holds for a
dimension $d-1$.
Let $B = [a_0,b_0]\times B'$ be a box in dimension $d$ [where $B'$ is
a $(d-1)$-dimensional box]. Let $\mu_0$ be the first marginal of $\mu$
(this is a measure on $[a_0,b_0]\subset\mathbb{R}$).\vspace*{1pt}

Let $b = \lfloor n^{1/d} \rfloor$ be the integer part of $n^{1/d}$,
and let
$b_0= 1/{{ (2^{1/d} - 1 )}}$. If \mbox{$b\leq b_0$}, we reason as
in the base case,
on the one-dimensional measure $\mu_0$, to find a partition of $B$ in $n$
slices of mass $\mu(B)/n$. Since the number of slices is less than the
constant $(b_0+1)^d$, the edge length is controlled as needed.

If\vspace*{1pt} $b>b_0$, we look for a decomposition of $n$ as a sum of $b$ integers
$n_i$, each as close to $n^{(d-1)/d}$ as possible: the idea is to cut $B$
along the first dimension in $b$ slices, and to apply the induction
hypothesis to cut the slice $i$ in $n_i$ parts.

To this end, decompose the integer $n$ in base $b$
\[
\exists\alpha_0, \alpha_1, \ldots,\alpha_{d}
\in\{0,\ldots, b-1\}^{d+1} \qquad n = \sum_{k=0}^{d}
\alpha_k b^k.
\]
The condition $b>b_0$ guarantees that $b+1 < 2^{1/d}b$, which implies
that $\alpha_d = 1$. Therefore,
\[
\exists\alpha_0, \alpha_1, \ldots,\alpha_{d-1}
\in\{0,\ldots, b-1\}^{d} \qquad n = b^d + \sum
_{k=0}^{d-1} \alpha_k b^k.
\]
Writing $\alpha_k = \sum_{i=1}^{b} \IND_{{ \{i \leq\alpha
_k \}}}$ we get
\begin{eqnarray*}
n &=& \sum_{i=1}^b {{
\Biggl(b^{d-1} + \sum_{k=0}^{d-1}
\IND_{{ \{
i\leq\alpha_k \}}} b^k \Biggr)}} = \sum
_{i=1}^b n_i,
\end{eqnarray*}
where $n_i = b^{d-1} + \sum_{k=0}^{d-1} \IND_{{ \{i\leq\alpha
_k \}}} b^k$.
From this expression, we get the bound
\(
b^{d-1}
\leq n_i
\leq(b^d - 1)/(b-1).
\)
Since\vspace*{2pt} $k-1\geq k/2$ whenever $k\geq2$, using the inequalities $b\leq n^{1/d}$
and $b \geq n^{1/d} - 1 \geq\frac{1}{2} n^{1/d}$, we get
\[
\frac{1}{2^{d-1}} n^{(d-1)/d} \leq n_i \leq2
n^{(d-1)/d}.
\]

Now let us cut $B$ along its first dimension.
Recall that $\mu_0$ is the first marginal of~$\mu$.
By continuity there exist quantiles $a_0 = q_0 <
q_1<\cdots< q_b = b_0$ such that
\[
\forall1\leq i\leq b \qquad\mu_1 \bigl([q_{i-1},q_i]
\bigr) = \mu\bigl([q_{i-1},q_i]\times B'
\bigr) = \frac{n_i}{n} \mu(B).
\]
We apply the induction hypothesis separately for each $1\leq i \leq b$,
to the \mbox{$(d-1)$-}di\-men\-sional box $B'$, with the measure
\begin{eqnarray*}
\mu_i(\cdot) &=& \mu\bigl([q_{i-1},q_i]
\times\cdot\bigr)
\end{eqnarray*}
and the integer $n_i$ to obtain a decomposition
$B' = \bigcup_{j=1}^{n_i} B'_{i,j}$ such that:
\begin{longlist}[(2)]
\item[(1)] the edge lengths $B'_{i,j}$ are controlled;
\item[(2)] $\mu_i(B'_{i,j}) = \frac{1}{n_i}\mu_i(B')$.
\end{longlist}
Finally, for
all $1\leq i \leq b$ and all $1\leq j\leq n_i$, let
\[
B_{i,j} = [q_{i-1},q_i] \times
B'_{i,j}.
\]

Let us check that the partition $B = \bigcup_i \bigcup_j B_{i,j}$
satisfies the requirements. By definition,
\[
\mu(B_{i,j}) = \mu_i \bigl(B'_{i,j}
\bigr) = \frac{1}{n_i} \mu_i \bigl(B' \bigr) =
\frac{1}{n}\mu(B),
\]
so the first requirement is met.
To control the edge lengths, first remark that
\begin{eqnarray*}
l(B) &=& \min\bigl( b_0 - a_0, l \bigl(B'
\bigr) \bigr),\qquad L(B) = \max\bigl( b_0 - a_0, L
\bigl(B' \bigr) \bigr),
\\
l(B_{i,j}) &=& \min\bigl(q_i - q_{i-1}, l
\bigl(B'_{i,j} \bigr) \bigr),\qquad L(B_{i,j}) = \max
\bigl(q_i - q_{i-1}, L \bigl(B'_{i,j}
\bigr) \bigr).
\end{eqnarray*}
By the induction hypothesis, the bounds on $n_i$ and the fact that
$L(B')\leq L(B)$ we get
\begin{eqnarray*}
L \bigl(B'_{i,j} \bigr) &\leq&\frac{C(d-1,\delta)}{n_i^{1/(d-1)}} L
\bigl(B' \bigr)
\\
&\leq&\frac{2 C(d-1,\delta)}{n^{1/d}} L(B).
\end{eqnarray*}
On the other hand, reasoning as in the proof of the base case,
\begin{eqnarray*}
(q_i - q_{i-1}) \bigl|B'\bigr| &\leq&
\frac{1}{\delta} \frac{n_i}{n}\mu(B),
\\
\mu(B) &\leq&\frac{1}{\delta} (b_0 - a_0)
\bigl|B'\bigr|
\end{eqnarray*}
so
\begin{eqnarray*}
(q_i - q_{i-1}) &\leq&(b_0 -
a_0) \delta^{-2} \frac{n_i}{n} \leq L(B)
\frac{2\delta^{-2}}{n^{1/d}}.
\end{eqnarray*}
Therefore $L(B_{i,j}) \leq C(d,\delta) n^{-1/d} L(B)$. The proof of the
lower bound on $l(B_{i,j})$ follows the same lines and is omitted. This
completes the induction step, and the lemma is proved.
\end{pf*}

%
%
\begin{corollary}[(Lower bound)]\label{corborneinf}
Under assumptions \textup{(H1)--(H3)}, if $\beta_N\gg N\log(N)$, then
for all $A\subset\mathcal{M}_1(\mathbb{R}^d)$, it holds
\[
\liminf_{N\to\infty}\frac{\log Z_NQ_N(A)}{\beta_N} \geq- \inf
\bigl\{
I(\eta);
\eta\in\operatorname{int}(A), \eta\ll\mathrm{Lebesgue} \bigr\}.
\]
\end{corollary}

\begin{pf}
Let $A\subset\mathcal{M}_1(\R^d)$ be a Borel set, and let $\eta\in
\operatorname{int}(A)$ be
absolutely continuous with respect to Lebesgue with density $h$ and such
that $I(\eta)<+\infty$. For some sequence ${(\varepsilon_n)}_{n\geq
1}$ converging
to $0$, let us define, for all $n\geq1$,
\[
\eta_n:= (1-\varepsilon_n)\nu_n +
\varepsilon_n \lambda_n,
\]
where $d\nu_n(x)=\frac{1}{C_n} \min(h(x); n) \IND_{[-n;n]^d}(x)
\,dx$ and
$d\lambda_n(x)=\frac{1}{(2n)^d} \IND_{[-n;n]^d}(x) \,dx$, where the
normalizing constant $C_n\to1$, when $n\to+\infty$.

According to point (3) of Lemma~\ref{lemgoodrate}, we see that
\[
I(\nu_n)<\infty,\qquad%
I(\lambda_n)<\infty,
\qquad%
\int\!\!\!\int\varphi(x,y) \,d\nu_n(x)\,d
\lambda_n(y)<\infty,
\]
where $\varphi(x,y):=\frac{1}{2}{{ (V(x)+V(y)+W(x,y) )}}$
(this function takes
its values in $(-\infty,+\infty]$ and is bounded from below thanks to
\textup{(H3)}; see the proof of Lemma~\ref{lemgoodrate}). It holds
\[
I(\eta_n) = (1-\varepsilon_n)^2 I(
\nu_n) %
+ 2\varepsilon_n(1-
\varepsilon_n)\int\!\!\!\int\varphi(x,y) \,d\nu_n(x)\,d
\lambda_n(y) %
+ \varepsilon_n^2I(
\lambda_n).
\]
Choose $\varepsilon_n$ converging to $0$ sufficiently fast so that the
last two
terms above converge to $0$ when $n \to\infty$. According to point
(1) of
Lemma~\ref{lemgoodrate}, $V\in\mathbf{L}^1(\mu)$ and $W\in
\mathbf{L}^1(\mu^2)$; it
follows then easily from the dominated convergence theorem that $I(\nu_n)
\to I(\eta)$ when $n\to\infty$ and that $\eta_n$ converges to $\eta
$ for
the weak topology.

Let $r>0$ be such that $B(\eta,2r)\subset A$; for all $n$ large enough,
$B(\eta_n,r)\subset B(\eta,2r)\subset A$. Since $\eta_n$ satisfies the
assumptions of Proposition~\ref{propborneinf}, we conclude that for $n$
large enough,
\[
\liminf_{N\to\infty} \frac{\log Z_NQ_N(A)}{\beta_N} %
\geq\liminf
_{N\to\infty} \frac{\log Z_NQ_N(B(\eta_n,r))}{\beta
_N} %
\geq- I(
\eta_n).
\]
Letting $n\to\infty$ and optimizing over $\{\eta\in A, \eta\ll
\mathrm{Lebesgue}\}$ gives the conclusion.
\end{pf}

\begin{pf*}{End of the proof of Theorem~\ref{thldp}}
The properties of $I_\star$ and the upper bound in point (2) are already
known. The lower bound of point (2) is given by Corollary~~\ref{corborneinf}.

To prove point (3), let $A\subset\mathcal{M}_1(\mathbb{R}^d)$ be
some Borel
set and take
$\mu\in\operatorname{int}(A)$. According to assumption \textup
{(H4)}, there
exists a
sequence of absolutely continuous probability measures $\nu_n$ converging
weakly to $\mu$ and such that $I(\nu_n) \to I(\mu)$, when $n\to
\infty
$. For
all $n$ large enough, $\nu_n \in A$ so applying Corollary~\ref{corborneinf}, we conclude that
\[
\liminf_{N\to\infty}\frac{\log Z_NQ_N(A)}{\beta_N}\geq-I(\nu_n).
\]
Letting $n\to\infty$ and then optimizing over $\mu\in\operatorname
{int}(A)$ we
arrive at
\[
\liminf_{N\to\infty}\frac{\log Z_NQ_N(A)}{\beta_N} %
\geq-\inf\bigl
\{I(\mu); \mu\in\operatorname{int}(A) \bigr\}.
\]
According to point (2) of Proposition~\ref{propreductionldp}, we conclude
that $Q_N$ obeys the full LDP.
\end{pf*}

\subsection{\texorpdfstring{Proof of the almost-sure convergence.}{Proof of the almost-sure convergence}}
Let us establish the last part of Theorem~\ref{thldp}. First note
that since
$I_\star$ has compact sublevel sets and is bounded from below,
$I_\star$
attains is infimum, so $I_{\min}$ is not empty. For an
arbitrary fixed real $\varepsilon>0$, consider the complement of the
$\varepsilon
$-neighborhood of
$I_{\min}$ for the Fortet--Mourier distance
\[
A_\varepsilon:=(I_{\min})_\varepsilon^c:= \bigl\{
\mu\in\mathcal{M}_1\dvtx{d_\mathrm{FM}}(\mu,I_{\min})>
\varepsilon\bigr\}.
\]
Since $I$ is lower semi-continuous, $c_\varepsilon:=\inf_{\mu\in
A_\varepsilon}I(\mu)>0$,
thus $\mathbb{P}(\mu_N\in A_\varepsilon)\leq\exp(-\beta_N
c_\varepsilon)$, by the upper
bound of the
full large deviation principle. By the first Borel--Cantelli lemma, it
follows that almost surely,
$\lim_{N\to\infty}{d_\mathrm{FM}}(\mu_N,\break I_{\min})=0$.

\subsection{\texorpdfstring{Sufficient conditions for \textup{(H4)}.}{Sufficient conditions for \textup{(H4)}}}

The following proposition gives several sufficient conditions under which
assumption \textup{(H4)} holds true. Even if some of these conditions are
quite general, it is an open problem to find an even more general and natural
condition. One may possibly find some inspiration in \cite{MR2858167}.

%
%
\begin{proposition}[{[Sufficient conditions for \textup{(H4)}]}]\label
{propassumptionH4}
Let $V\dvtx\R^d\to\R$ and $W\dvtx\mathbb{R}^d\times\mathbb
{R}^d\to(-\infty,+\infty]$ be
symmetric, finite on $\mathbb{R}^d\times\mathbb{R}^d\setminus\{
(x,x); x\in\R
^d\}$ and
such that \textup{(H2)} and \textup{(H3)} hold true. Assumption
\textup{(H4)} holds in each of the following cases:
\begin{longlist}[(3)]
\item[(1)] $W$ is finite and continuous on $\R^d\times\R^d$;
\item[(2)] for all $x\in\mathbb{R}^d$, the function $y\mapsto W(x,y)$ is super
harmonic, that is, $W$~satisfies
\[
W(x,y)\geq\frac{1}{|B(y,r)|}\int_{B(y,r)} W(x,z) \,dz\qquad\forall r>0,
\]
where $|B(y,r)|$ denotes the Lebesgue measure of the ball of center $y$
and radius~$r$;
\item[(3)] the function $W$ is such that $W(x+a,y+a)=W(x,y)$ for all $x,y,a
\in\mathbb{R}^d$ and the function $J$ defined by
%
%
\begin{equation}
\label{eqdefJ} J(\mu)=\int\!\!\!\int W(x,y) \,d\mu(x)\,d\mu(y)
\end{equation}
is convex on the set of compactly supported probability measures.
\end{longlist}
\end{proposition}

\begin{pf}
Let $\mu\in\mathcal{M}_1(\mathbb{R}^d)$ be such that $I(\mu
)<\infty$.
Recall that,
according to point (1) of Lemma~\ref{lemgoodrate}, under assumptions
\textup{(H2)--(H3)}, the condition $I(\mu)<\infty$ implies that
\[
\int|V| \,d\mu<+\infty\quad\mbox{and}\quad\int\!\!\!\int|W| \,
d\mu^2<+
\infty.
\]
Moreover, it follows from \textup{(H2)} and \textup{(H3)}
that $W$
is bounded from below on every compact, and so the definition (\ref
{eqdefJ}) of $J(\mu)$ makes sense if $\mu$ is compactly supported.

For all $R>0$, let us define $\mu_R$ as the normalized restriction of
$\mu$
to $[-R;R]^d$. Using the dominated convergence theorem and point (1) of
Lem\-ma~\ref{lemgoodrate}, it is not hard to see that $\mu_R$ converges
weakly to $\mu$ and that $I(\mu_R)\to I(\mu)$ when \mbox{$R \to+\infty$}. To
regularize $\mu_R$, we consider $\mu_{R,\varepsilon} =\operatorname{Law}
(X_R+\varepsilon U)$,
$\varepsilon\leq1$, where $X_R$ is distributed according to $\mu_R$ and
$U$ is
uniformly distributed on the Euclidean unit ball $B_1$ of $\mathbb
{R}^d$. It is
clear that $\mu_{R,\varepsilon}$ has a density with respect to
Lebesgue measure.
Moreover, $\mu_{R,\varepsilon}\to\mu_R$, when $\varepsilon\to0$.
Indeed, if $f\dvtx\mathbb{R}^d
\to\R$ is continuous, it is bounded on $[-R;R]^d+B_1$, and it follows that
\[
\int f \,d\mu_{R,\varepsilon}=\mathbb{E} \bigl[f(X_R+\varepsilon U)
\bigr]\mathop{\to}_{\varepsilon\to0}\mathbb{E} \bigl[f(X_R)
\bigr].
\]
This applies in particular to $f=V$. Now let us show in each cases that
$J(\mu_{R,\varepsilon})$ converges to $J(\mu_{R})$, when
$\varepsilon$ goes to
$0$. Let
us write
\[
J(\mu_{R,\varepsilon}) = \mathbb{E} \bigl[W(X_R+\varepsilon
U,Y_R+\varepsilon V) \bigr],
\]
where $X_R,Y_R, U, V$ are independent and such that $Y_R\stackrel{(d)}{=}
X_R$ and $V \stackrel{(d)}{=}U$.
\begin{longlist}[(3)]
\item[(1)] If $W$ is finite and continuous on $\mathbb{R}^d\times\mathbb
{R}^d$, then
using the
boundedness of $W$ on $([-R;R]^d+B_1) \times([-R;R]^d+B_1)$, it follows
that $J(\mu_{R,\varepsilon})\to J(\mu_R)$ when $\varepsilon\to0$.
\item[(2)] If $W$ is super harmonic, then $W_\varepsilon(x,y):=\mathbb
{E}_{U,V}[W(x+\varepsilon U,y
+ \varepsilon V)] \leq W(x,y)$ for all $x,y$. Moreover, it follows from the
continuity of $W$ outside the diagonal that, for all $x\neq y$,
$W_\varepsilon(x,y) \to W(x,y)$ when $\varepsilon\to0$. Since $I(\mu
)<+\infty$,
$\mu$ does not have atoms, and so the diagonal is of measure $0$ for
$\mu^2$. It follows from the dominated convergence theorem that
$J(\mu_{R,\varepsilon}) \to J(\mu_R)$ as $\varepsilon\to0$.
\item[(3)] Denoting by $\mu_R^x$ the law of $X_R+x$, we see that $\mu
_{R,\varepsilon} =
\mathbb{E}_U[\mu_R^{\varepsilon U}]$. Therefore, the convexity of $J$
yields to
\begin{eqnarray*}
J(\mu_{R,\varepsilon}) %
&\leq& \mathbb{E}_U \bigl[ J \bigl(
\mu_R^{\varepsilon U} \bigr) \bigr] %
\\
&=&
\mathbb{E}_U \biggl[\int\!\!\!\int W(x+\varepsilon U,y+\varepsilon U)
\,d\mu_R(x)\,d \mu_R(y) \biggr] %
= J(
\mu_R),
\end{eqnarray*}
where the last equality comes from the property $W(x+a,y+a)=W(x,y)$. On
the other hand, Fatou's lemma implies that $\liminf_{\varepsilon\to0}
J(\mu_{R,\varepsilon}) \geq J(\mu_R)$. Therefore $J(\mu
_{R,\varepsilon}) \to
J(\mu_R)$, when $\varepsilon$ goes to $0$.
\end{longlist}
We conclude from the above discussion that for any $\delta>0$, it is
possible to choose $R$ sufficiently large and then $\varepsilon$ sufficiently
small so that ${d_\mathrm{FM}}(\mu_{R,\varepsilon},\mu)\leq\delta
$ and
$|I(\mu_{R,\varepsilon
}) -
I(\mu)| \leq\delta$. This completes the proof.
\end{pf}

\section{\texorpdfstring{Tools from potential theory.}{Tools from potential theory}}\label{secpotential}
In this section, we recall results from potential theory that will prove
useful when we discuss the proof of Theorem~\ref{thriesz} and
Corollary~\ref{cocoulomb}. There are many textbooks on potential
theory, with
different point of views; our main source is \cite{MR0350027}, where
the Riesz
case is well-developed.

In this
section, and unless otherwise stated, we set $k_\alpha:=k_{\Delta
_\alpha}$ and
we take
$W(x,y):=k_\alpha(x-y)$, $0<\alpha<d$, $d\geq1$. We denote,
respectively, by
$\mathcal{M}_1\subset\mathcal{M}_\infty\subset\mathcal
{M}_+\subset\mathcal{M}_{\pm}$ the sets of
probability measures, of positive measures integrating
$k_\alpha(\cdot)\IND_{{{\vert\cdot\vert}}>1}$, of positive
measures, and of signed
measures on $\mathbb{R}^d$.

\subsection{\texorpdfstring{Potentials and interaction energy.}{Potentials and interaction energy}}
We benefit from the constant sign of the Riesz kernel: $k_\alpha\geq
0$, contrary
to the Coulomb kernel in dimension $d=2$ and its logarithm. Following
\cite{MR0350027}, page~58, the \emph{potential} of $\mu\in\mathcal
{M}_+$ is the function
$U^\mu_\alpha\dvtx\mathbb{R}^d\to[0,\infty]$ defined for every
$x\in\mathbb{R}
^d$ by
%
%
\begin{equation}
\label{eqpotU} U^\mu_\alpha(x):=\int W(x,y) \,d\mu(y)
=\int k_\alpha(x-y) \,d\mu(y).
\end{equation}
Note that $U^\mu_\alpha(x)=\infty$ if $\mu$ has a Dirac mass at
point $x$. By
using the Fubini theorem, for every $\mu\in\mathcal{M}_+$, we have
$U^\mu_\alpha
<\infty$
Lebesgue almost everywhere if and only if $\mu\in\mathcal{M}_\infty
$. This
explains actually the
condition $0<\alpha<d$ taken in the Riesz potential, which is related
to polar
coordinates ($dx=r^{d-1} \,dr\,d\sigma_{ d}$). In fact if $\mu\in
\mathcal{M}
_\infty$, then
$U^\mu_\alpha$ is a locally Lebesgue integrable function. Moreover,
as Schwartz
distributions, we have $U^\mu_\alpha=k_\alpha*\mu$ and, with the
notation of
(\ref{eqfs-riesz}),
\[
-c_\alpha\Delta_\alpha U^\mu_\alpha=(-c_\alpha
\Delta_\alpha k_\alpha)*\mu=\mu.
\]
The interaction energy is the quadratic functional
$J_\alpha\dvtx\mathcal{M}_+\mapsto[0,\infty]$ defined by
\[
J_\alpha(\mu):=\int\!\!\!\int W(x,y) \,d\mu(x)\,d\mu(y) =\int
U^\mu_\alpha \,d\mu.
\]
Note that $J_\alpha(\mu)=\infty$ if $\mu$ has a Dirac mass, and in
particular
$J_\alpha(\mu_N)=\infty$.

In the Coulomb case where $\alpha=2$, we have $c_2J_2(\mu)=-\int
U_2^\mu
\Delta
U_2^\mu \,dx=\int{{\vert\nabla U_2^\mu\vert}}^2 \,dx$. The
quantity $\nabla
U_2^\mu$
is the (electric) field generated by the (Coulomb) potential $U_2^\mu
$, and
this explains the term ``carr\'e-du-champ''\break  (``square of the field'' in
French) used for $J_2(\mu)$.

%
%
\begin{lemma}[(Positivity and convexity on $\mathcal{M}_+$)]\label
{leconvexity} %
\begin{itemize}
\item for every $\mu\in\mathcal{M}_+$ we have $J_\alpha(\mu)\geq
0$ with equality
if and only if
$\mu=0$;
\item$J_\alpha\dvtx\mathcal{M}_+\mapsto[0,\infty]$ is strictly convex:
for every
$\mu,\nu\in\mathcal{M}_+$ with $\mu\neq\nu$, we have
\[
\forall t\in(0,1)\qquad J_\alpha\bigl(t\mu+(1-t)\nu\bigr)<
tJ_\alpha(\mu)+(1-t)J_\alpha(\nu);
\]
\item$\mathcal{E}_{\alpha,+}:=\{\mu\in\mathcal{M}_+\dvtx J_\alpha
(\mu
)<\infty\}$ is a convex cone.
\end{itemize}
\end{lemma}

We recall that in classical harmonic analysis, a function
$K\dvtx\mathbb{R}\times\mathbb{R}\to\mathbb{R}$ is called a \emph
{positive definite
kernel} when
$\sum_{i=1}^nx_iK(x_i,x_j)\bar{x}_j\geq0$ for every $n\geq1$ and every
$x\in\mathbb{C}^n$.
If this holds only when $x_1+\cdots+x_n=0$, the kernel is said to be
\emph{weakly positive definite}. The famous Bochner theorem states
that a
kernel is positive definite if and only if it is the Fourier transform
of a
\emph{finite} Borel measure. The famous Schoenberg theorem states for
every $f\dvtx\mathbb{R}_+\to\mathbb{R}_+$, the kernel $(x,y)\mapsto
f({{\vert x-y
\vert}}^2)$ is positive
definite on $\mathbb{R}^d$ for every $d\geq1$ if and only if $f$ is
the Laplace
transform of a \emph{finite} Borel measure on $\mathbb{R}_+$. The famous
Bernstein theorem states that if $f\dvtx\mathbb{R}\to\mathbb{R}$ is
continuous and
$\mathcal{C}^\infty((0,\infty))$, then $f$ is the Laplace transform
of a finite Borel
measure on $\mathbb{R}_+$ if and only if $f$ is \emph{completely monotone}:
$(-1)^nf^{(n)}\geq0$ for every $n\geq0$. For all these notions, we
refer to
\cite{MR747302,MR2132704}.

The proof of Lemma~\ref{leconvexity} is short and self-contained. It relies
on the fact that the convexity of the functional is equivalent to the fact
that $W$ is a weakly positive definite kernel, which is typically the case
when $W$ is a mixture of shifted Gaussian kernels, which are the most useful
weakly positive definite kernels. For example, this works if for some arbitrary
measurable $\alpha,\beta\dvtx\mathbb{R}\to\mathbb{R}$ and Borel
measure $\eta
$,\vadjust{\goodbreak} and every
$x,y\in\mathbb{R}^d$,
\[
W(x,y)=w \bigl({{\vert x-y \vert}} \bigr)=\int_0^\infty{{
\bigl(e^{-\alpha
^2(t){{\vert x-y \vert}}^2}+\beta(t) \bigr)}} \,d\eta(t).
\]
The shift $\beta$ can be $<0$, which allows nonpositive definite
kernels such
as the logarithmic kernel (note that the Riesz kernel is positive definite).
The method is used for the logarithmic kernel in \cite{MR1465640},
Proof of Property~2.1(4), with the following mixture:
\[
\log\frac{1}{{{\vert x-y \vert}}}%
=\int_0^\infty
\frac{1}{2t}{{ \bigl(e^{-{{{\vert x-y \vert
}}^2}/{(2t)}}-e^{-{1}/{(2t)}} \bigr)}} \,dt.
\]
This kernel has a sign change and a double singularity near zero and infinity,
which can be circumvented by using a cutoff. Alternatively, one may
proceed by
regularization and use the Bernstein theorem with the completely monotone
function $f(t)=(\varepsilon+t)^{-\beta}$, $\beta,\varepsilon>0$,
and then the Schoenberg
theorem; see, for example,~\cite{MR1422615}. For instance, for the
logarithmic kernel,
the following representation is used in \cite{MR1746976}, Chapter~5:
\[
\log\frac{1}{\varepsilon+{{\vert x-y \vert}}} =\int_0^\infty{{
\biggl(
\frac{1}{\varepsilon+1+{{\vert x-y
\vert}}}-\frac{1}{1+t} \biggr)}} \,dt.
\]
Finally, let us mention that for the Riesz kernel, yet another short
proof of
Lemma~\ref{leconvexity}, based on the formula $k_\alpha=ck_{\alpha
/2}*k_{\alpha/2}$,
can be found in \cite{MR0350027}, Theorem~1.15, page~79.

\begin{pf*}{Proof of Lemma~\ref{leconvexity}}
Set $\beta:=d-\alpha$. We start from the identity
\[
\Gamma(1+\alpha)=c^{1+\alpha}\int_0^\infty
t^\alpha e^{-ct} \,dt,\qquad c>0, \alpha>-1.
\]
Taking $c={{\vert x-y \vert}}^2$ and $1+\alpha=\beta/2$, we get,
for every $x,y\in\mathbb{R}^d$,
\[
k_\alpha(x-y)%
= \int_0^\infty
f(t) e^{-t{{\vert x-y \vert}}^2} \,dt \qquad\mbox{where }
f(t):=\frac{t^{\beta/2-1}}{\Gamma(\beta/2)}.
\]
Now for every $\mu\in\mathcal{M}_+$ such that $J_\alpha(\mu
)<\infty$,
\[
J_\alpha(\mu) =\int_0^\infty f(t){{
\biggl(\int\!\!\!\int e^{-t{{\vert x-y \vert}}^2} \,d\mu(x)\,d\mu(y)
\biggr)}} \,dt.
\]
Expressing the Gaussian kernel as the Fourier transform of a Gaussian
kernel, we get, by writing $e^{i{{ \langle x-y,w \rangle
}}}=e^{i{{ \langle x,w \rangle}}}e^{-i{{ \langle
y,w \rangle}}}$
and using the Fubini theorem,
\begin{eqnarray*}
&& \int\!\!\!\int e^{-t{{\vert x-y \vert}}^2} \,d\mu(x)\,d\mu(y)
\\
&&\qquad = (4\pi
t)^{-d/2} \int\!\!\!\int%
{{ \biggl(\int_{\mathbb{R}^d}
e^{i{{ \langle x-y,w \rangle
}}}e^{(1/(4t)){{\vert w
\vert}}^2} \,dw \biggr)}} \,d\mu(x)\,d\mu(y)
\\
&&\qquad = (4\pi t)^{-d/2} %
\int_{\mathbb{R}^d}
\underbrace{{{ \biggl\vert\int e^{i{{ \langle
x,w \rangle}}} \,d\mu(x) \biggr\vert
}}^2}_{K_w(\mu)} e^{-(1/(4t)){{\vert w \vert}}^2} \,dw.
\end{eqnarray*}
Now $K_w$ is clearly convex since for every $\mu,\nu\in\mathcal
{M}_1(\mathbb{R}^d)$
and every
$t\in(0,1)$,
\begin{eqnarray*}
&& \frac{tK_w(\mu)+(1-t)K_w(\nu)-K_w(t\mu+(1-t)\nu)}{t(1-t)}
\\
&&\qquad =K_w(\mu
-\nu)
={{ \biggl\vert\int
e^{i{{ \langle x,w \rangle}}} \,d(\mu-\nu) (x) \biggr\vert}}^2\geq0.
\end{eqnarray*}
It follows then that $J_\alpha$ is also convex as a conic combination
of convex
function. Let us establish now the strict convexity of $J_\alpha$. Let us
suppose that $\mu,\nu\in\mathcal{M}_1(\mathbb{R}^d)$ with
$J_\alpha(\mu
)<\infty$ and
$J_\alpha(\nu)<\infty$ and $tJ_\alpha(\mu)+(1-t)J_\alpha(\nu
)=J_\alpha(t\mu
+(1-t)\nu)$
for some $t\in(0,1)$. Then
\[
J_\alpha(\mu-\nu) = \frac{tJ_\alpha(\mu)+(1-t)J_\alpha(\nu
)-J_\alpha
(t\mu+(1-t)\nu
)}{t(1-t)} =0.
\]
Arguing as before, we find
\[
0= J_\alpha(\mu-\nu) = \int_0^\infty
f(t){{ \biggl[(4\pi t)^{-d/2} \int_{\mathbb{R}^d}
K_w(\mu-\nu) e^{-({1}/(4t)){{\vert w \vert}}^2} \,dw \biggr]}} \,dt.
\]
Hence, for every $t>0$ (a single $t>0$ suffices in what follows),
\[
\int_{\mathbb{R}^d} K_w(\mu-\nu)
e^{-({1}/(4t)){{\vert w \vert}}^2} \,dw=0.
\]
Thus, the Fourier transform of $\mu-\nu$ vanishes almost everywhere, and
therefore $\mu=\nu$.

Finally, $\mathcal{E}_{\alpha,+}$ is clearly a cone, and its
convexity comes from the
convexity of~$J_\alpha$.
\end{pf*}

Following\vspace*{1pt} \cite{MR0350027}, page~62, for any $\mu=\mu_+-\mu_-\in
\mathcal{M}_{\pm
}$ such
that $\mu_\pm\in\mathcal{M}_\infty$, we have $U^{\mu_\pm}_\alpha
<\infty$ Lebesgue almost
everywhere, and we may define for Lebesgue almost every $x$
\[
U^\mu_\alpha(x):=U^{\mu_+}_\alpha(x)-U^{\mu_-}_\alpha(x)
\in(-\infty,+\infty).
\]
Following \cite{MR0350027}, page~77, for every $\mu=\mu_+-\mu_-\in
\mathcal{M}
_{\pm}$
such that $\mu_\pm\in\mathcal{M}_\infty$ and
\[
\int U^{\mu_+}_\alpha \,d\mu_-<\infty\quad\mbox{and}\quad\int
U^{\mu_-}_\alpha \,d\mu_+<\infty,
\]
we may define $J_\alpha(\mu)\in(-\infty,+\infty]$ as (thanks to
the Fubini
theorem)
\[
J_\alpha(\mu):=\int U^\mu_\alpha \,d\mu=\int
U^{\mu_+}_\alpha \,d\mu_++\int U^{\mu_-}_\alpha d
\mu_- -\int U^{\mu_+}_\alpha \,d\mu_--\int U^{\mu_-}_\alpha
\,d\mu_+.
\]
More generally, for every $\mu_1,\mu_2\in\mathcal{M}_{\pm}$ such that
$\mu_{1\pm},\mu_{2\pm}\in\mathcal{M}_\infty$ and
\[
\int U^{\mu_{1\pm}}_\alpha \,d\mu_{2\mp}<\infty,
\]
we may define $J_\alpha(\mu_1,\mu_2)\in(-\infty,+\infty]$ by
\begin{eqnarray*}
J_\alpha(\mu_1,\mu_2)&:=&\int
U^{\mu_1}_\alpha \,d\mu_2
\\
&=& \int U^{\mu_{1+}}_\alpha
\,d\mu_{2+}+\int U^{\mu_{1-}}_\alpha \,d\mu_{2-} -
\int U^{\mu_{1+}}_\alpha \,d\mu_{2-}-\int
U^{\mu_{1-}}_\alpha \,d\mu_{2+}.
\end{eqnarray*}
Following \cite{MR0350027}, page~78, since $k_\alpha$ is symmetric,
then the
\emph{reciprocity law} holds,
\[
J_\alpha(\mu_1,\mu_2)=J_\alpha(
\mu_2,\mu_1), \qquad\mbox{i.e., } \int
U^{\mu_1}_\alpha \,d\mu_2=\int U^{\mu_2}_\alpha
\,d\mu_1.
\]
Let $\mathcal{E}_\alpha$ be the set of elements of $\mathcal{M}_{\pm
}$ for which $J_\alpha$ makes
sense and is finite. As pointed out by Landkof \cite{MR0350027} in his
preface, a very nice idea going back to Cartan consists of seeing
$J_\alpha$
as a Hilbert structure on $\mathcal{E}_\alpha$. This idea is simply
captured by the
following lemma, which is the analogue of Lemma~\ref{leconvexity} for signed
measures of finite energy.

%
%
\begin{lemma}[{[Properties of $(\mathcal{E}_\alpha,J_\alpha)$]}]\label
{lehilbert}
%
\begin{itemize}
\item$J_\alpha$ is lower semi-continuous on $\mathcal{E}_\alpha$
for the vague topology
(i.e., with respect to continuous functions with compact support);
\item$\mathcal{E}_\alpha$ is a vector space and $(\mu_1,\mu
_2)\mapsto
J_\alpha(\mu_1,\mu_2)$ defines a scalar product on $\mathcal
{E}_\alpha$.
\end{itemize}
In particular for every $\mu\in\mathcal{E}_\alpha$, we have
$J_\alpha(\mu)=J_\alpha(\mu,\mu)\geq0$ with equality if and only
if $\mu=0$;
and moreover,
$J_\alpha\dvtx\mathcal{E}_\alpha\mapsto(-\infty,\infty)$ is strictly
convex: for every
$\mu,\nu\in\mathcal{E}_\alpha$ with $\mu\neq\nu$,
\[
\forall t\in(0,1)\qquad\frac{tJ_\alpha(\mu)+(1-t)J_\alpha(\nu
)-J_\alpha(t\mu+(1-t)\nu)}{t(1-t)} =J_\alpha(\mu-\nu)>0.
\]
\end{lemma}

\begin{pf}
The lower semi-continuity for the vague convergence follows from the fact
that $k_\alpha\geq0$; see \cite{MR0350027}, page~78. The vector space
nature of
$\mathcal{E}_\alpha$ is immediate from its definition. The
bilinearity of
$(\mu_1,\mu_2)\mapsto J_\alpha(\mu_1,\mu_2)$ is immediate. By reasoning
as in
the proof of Lemma~\ref{leconvexity}, we get $J_\alpha(\mu,\mu
)\geq0$ for
every $\mu\in\mathcal{E}_\alpha$, with equality if and only if $\mu=0$.
\end{pf}

Following \cite{MR0350027}, Theorems 1.18~and~1.19, page~90, for this
pre-Hilbertian topology, it can be shown that $\mathcal{E}_{\alpha,+}$
is complete while
$\mathcal{E}_\alpha$ is not complete if $\alpha>1$, and that
$J_\alpha$ is not
continuous for
the vague topology.

\subsection{\texorpdfstring{Capacity and ``approximately/quasi-everywhere.''}{Capacity and ``approximately/quasi-everywhere''}}
The notion of capacity is central in Potential Theory. We just
need basic facts on zero-capacity sets. Once more we follow
the presentation of Landkof \cite{MR0350027}, Chapter II.1,
to which we refer for additional details, references and proofs.

For any compact set $K$, consider the minimization problem
\[
W_\alpha(K) = \inf\bigl\{ J_\alpha(\nu); \nu\in
\mathcal{M}_1 \cap\mathcal{E}_\alpha, \operatorname{supp}(
\nu)\subset K \bigr\}.
\]
The boundedness of $K$ implies that $W_\alpha(K)\in(0,\infty]$. Its
inverse $C_\alpha(K)$ is
called the \emph{capacity} of the compact set $K$. The capacity of $K$
is zero if and only if
there is no measure of finite energy supported in $K$.

On general sets on can define an ``inner capacity'' and an ``outer
capacity'' by
\begin{eqnarray*}
\underbar{C}_\alpha(A) &=& \sup\bigl\{ C_\alpha(K), K\subset A,
K \mbox{ compact} \bigr\},
\\
\widebar{C}_\alpha(A) &=& \inf\bigl\{ \underbar{C}_\alpha(O), A
\subset O, O \mbox{ open} \bigr\}.
\end{eqnarray*}
It can be shown (see \cite{MR0350027}, Theorem~2.8) that if $A$ is a
Borel set,
these two quantities coincide---$A$ is said to be ``capacitable'' and the
common value is called the capacity of $A$.

A property $P(x)$ is said to hold ``approximately everywhere'' if the
set $A$
of $x$ such that $P(x)$ is false, has zero inner capacity, and
``quasi-everywhere'' if it has zero outer capacity. For many ``reasonable''
$P(x)$, the set $A$ is Borel, and the two notions coincide. The following
result \cite{MR0350027}, Theorems 2.1~and~2.2, shows that, for such
``reasonable'' properties, ``quasi-everywhere'' means ``$\nu$-almost surely,
for all measures $\nu$ of finite energy.''

%
%
\begin{theorem}[(Zero capacity Borel sets)]\label{propqe}
A Borel set $A$ has zero capacity if and only if, for any measure $\nu
$ of
finite energy, $\nu(A) = 0$. In particular, if $C_\alpha(A) > 0$, $A$
has a
positive inner capacity, and there exists a compact $K\subset A$ and a
probability measure $\nu$ of finite energy such that $\operatorname
{supp}(\nu)
\subset K$.
\end{theorem}

\subsection{\texorpdfstring{The Gauss averaging principle.}{The Gauss averaging principle}}
In the classical Coulombian case ($\alpha= 2$), we will need
the following result, known as Gauss's averaging
principle. In $\mathbb{R}^d$, for all $r>0$, let $\sigma_r$ be the
surface measure
on the sphere $\partial B(0,r)$; its total mass is $\sigma_d r^{d-1}$
where $\sigma_d$ is the surface of the unit sphere.\vadjust{\goodbreak}

%
%
\begin{theorem}[(Gauss's averaging principle)]\label{thmgauss}
In $\mathbb{R}^d$,
\[
\frac{1}{r^{d-1}\sigma_d}%
\int_{\partial B(0,r)} \frac{1}{{{\vert x-y \vert}}^{d-2}}
\,d\sigma_r(y) %
= %
\cases{ \displaystyle
\frac{1}{r^{d-2}}, &\quad if ${{\vert x \vert}} < r$,
\vspace*{3pt}\cr
\displaystyle
\frac{1}{{{\vert x \vert}}^{d-2}}, &\quad if ${{\vert x \vert
}} > r$.} %
\]
\end{theorem}

This result can be found in \cite{Hel09}, Lemma~1.6.1, page 21.

\section{\texorpdfstring{Proof of the properties of the minimizing measure.}{Proof of the properties of the minimizing measure}}\label{secminimizing}

The proof of Theorem~\ref{thriesz} is decomposed in two steps. We
begin by
proving the existence, uniqueness and the support properties of $\mu
_\star$ in
Section~\ref{secexistenceunicite}. The characterization of $\mu
_\star
$ is
proved in Section~\ref{seccarac}.

Recall that, for a probability measure $\mu$, we have defined
\[
I(\mu) = \frac{1}{2}J_\alpha(\mu) + \int V \,d\mu.
\]
In this section we consider the following minimization problem:
%
%
\begin{equation}
\label{eqproblemeOrigine} \mathcal{P}\dvtx\inf{{ \bigl\{I(\mu),
\mu\in
\mathcal{M}_1
\bigr\}}}.
\end{equation}

\subsection{\texorpdfstring{Existence, uniqueness and compactness of the support.}{Existence, uniqueness and compactness of the support}}\label{secexistenceunicite}
The existence of a minimizer for $\mathcal{P}$ is clear since we have
already seen
that $I$ has compact level sets.

Since $I(\mu)< \infty$ implies that $\mu\in\mathcal{E}_\alpha$
and $\int V
\,d\mu< \infty$,
the problem $\mathcal{P}$ is equivalent to
%
%
\begin{equation}
\label{eqproblemeRestreint} \mathcal{P}_\alpha\dvtx\inf{{ \bigl\{
I(\mu),
\mu\in\mathcal
{M}_1\cap\mathcal{E}_\alpha\mbox{ such that } V \in
L^1(\mu) \bigr\}}}
\end{equation}
in that they have the same values and the same minimizers.
Let us call $p$ the common value.

Suppose $\mu$ and $\nu$ are two measures in $\mathcal{M}_1\cap
\mathcal{E}_\alpha$ such that
$V\in\rL^1(\mu)\cap\rL^1(\nu)$. Let $\psi\dvtx t\in[0,1]\mapsto
[0,\infty)$ by
%
%
\begin{eqnarray}
\label{eqdefPsi} \psi(t)&:=&I \bigl((1-t)\mu+t\nu\bigr)
\nonumber\\[-8pt]\\[-8pt]
&=&\frac{1}{2} J_\alpha\bigl( (1-t)\mu+ t \nu\bigr) +
(1-t)\int V \,d\mu+ t\int V \,d\nu.\nonumber
\end{eqnarray}
By Lemma~\ref{lehilbert}, $\psi$ is strictly convex if $\mu\neq\nu$.
If $\mu$ and $\nu$ minimize $I$, then they are in $\mathcal{M}_1\cap
\mathcal{E}_\alpha$, so
$\psi$ is well defined, and since $\psi(0) = I(\mu) = I(\nu)= \psi(1)$,
$\mu$ must be equal to $\nu$. Therefore the minimizer $\mu_\star$
is unique.

Let us now prove that $\mu_\star$ has compact support.
This result also holds in dimension $2$ with the logarithmic potential;
see \cite{MR1485778}, Theorem~1.3, page~27.
To this end, let us define, for any compact $K$, a new minimization problem,
\[
\mathcal{P}_K\dvtx \inf{{ \bigl\{I(\mu), \mu\in\mathcal{M}_1
\cap\mathcal{E}_\alpha, \operatorname{supp}(\mu) \subset K \bigr\}}},
\]
and let $p_K$ be the value of $\mathcal{P}_K$.

%
%
\begin{lemma}[(Reduction to restricted optimization problem)]
Let $K$ be a compact set, and suppose that $V(x) \geq2 p + 3$ when
$x\notin
K$, where $p$ is the common value of $\mathcal{P}$, $\mathcal
{P}_\alpha$ [defined by
(\ref{eqproblemeOrigine}) and (\ref{eqproblemeRestreint})]. Then the
problems $\mathcal{P}$ and $\mathcal{P}_K$ are equivalent: their
values $p$ and $p_K$ are
equal, the minimizer exists and is the same. In particular, the minimizer
$\mu_\star$ of the original problem $\mathcal{P}$ satisfies
$\operatorname{supp}(\mu_\star)
\subset K$.
\label{lecompactSupport}
\end{lemma}
\begin{pf}
Suppose $\mu$ is such that $I(\mu) \leq p + 1$. We will prove that, if
$\mu(K) < 1$, we can find a $\mu_K$, supported in $K$ such that
$I(\mu
_K) <
I(\mu)$. This clearly implies that the two values $p_K$ and $p$ coincide.
Since we know that the minimizer $\mu_\star$ of the original problem exists,
this also proves that it must be supported in $K$.

Let us now construct $\mu_K$ as the renormalized restriction of $\mu$ to
$K$. First, remark that $\mu(K)$ cannot be zero, since
\[
p + 1 \geq I(\mu) \geq\bigl(1- \mu(K) \bigr) (2 p + 3).
\]
Therefore we can define
\[
\mu_K(A) = \frac{1}{\mu(K)} \mu(K\cap A).
\]
Since by assumption $\mu(K) < 1$, we may similarly define $\mu
_{K^c}$. The
measure $\mu$ is the convex combination
\[
\mu= \mu(K) \mu_K + \bigl(1 - \mu(K) \bigr) \mu_{K^c}.
\]

The positivity of $V$, $W$ and the choice of $K$ imply that
\begin{eqnarray*}
I(\mu) &=& \frac{1}{2} J_\alpha(\mu) %
+ \mu(K) \int V \,d
\mu_K + \bigl(1-\mu(K) \bigr)\int V \,d\mu_{K^c}
\\
&\geq&\frac{1}{2} \mu(K)^2 J_\alpha(\mu_K)
+ \mu(K)^2 \int V \,d\mu_K + \bigl(1-\mu(K)
\bigr) (2p + 3),
\end{eqnarray*}
since $J_\alpha(\mu_{K^c})$ and the interaction energy
$J_\alpha(\mu_K,\mu_{K^c})$ are both nonnegative. Therefore
\[
I(\mu) \geq\mu(K)^2 I(\mu_K) + \bigl(1 - \mu(K) \bigr)
(2p + 3).
\]
Assume that $I(\mu_K)\geq I(\mu)$. Then
\[
I(\mu) \bigl(1-\mu(K)^2 \bigr) \geq\bigl(1 - \mu(K) \bigr) (2p +
3).
\]
Using the fact that $I(\mu) \leq p + 1$, and dividing by $1 - \mu
(K)$, we
get
\[
2(p + 1) \geq(p + 1) \bigl(1 + \mu(K) \bigr) \geq2 p + 3,
\]
a contradiction. Therefore $I(\mu_K) < I(\mu)$, and the proof is complete.
\end{pf}

\subsection{\texorpdfstring{A criterion of optimality.}{A criterion of optimality}}\label{seccarac}

In this section we prove the items (5), (6)~and~(7) of Theorem~\ref{thriesz}. The corresponding result in dimension $2$ for the logarithmic
potential can be found in \cite{MR1485778}, Theorem~3.3, page~44. %
We adapt it, using fully the pre-Hilbertian structure rather than the
principle of domination when it is possible.

\begin{pf*}{Proof of item (5) of Theorem~\ref{thriesz}}
We already know that $\mu_\star$ has compact support. The first step
is to
show that $\mu_\star$ satisfies (\ref{eqrob1}) and (\ref{eqrob2}). Let
$\mu= \mu_\star$, and let $\nu$ be in $\nu\in\mathcal{M}_1\cap
\mathcal{E}_\alpha$ such that
$V\in\rL^1(\nu)$. Recall the function $\psi$ from (\ref{eqdefPsi}),
\[
\psi(t) = I \bigl(( 1-t)\mu_\star+ t\nu\bigr).
\]
Since $J_\alpha$ is quadratic we get
\begin{eqnarray*}
\psi(t) &=&\int V \,d\mu_\star+ t\int V \,d(\nu-\mu_\star) +
\frac{1}{2}J_\alpha\bigl(\mu_\star+ t(\nu-
\mu_\star) \bigr)
\\
&=& \int V \,d\mu_\star+ t\int V \,d(\nu-\mu_\star)
\\
&&{} +
\frac{1}{2} \bigl(J_\alpha(\mu_\star) + t^2
J_\alpha(\nu-\mu_\star)
+ 2t J_\alpha(
\mu_\star,\nu-\mu_\star) \bigr).
\end{eqnarray*}
Therefore,
%
%
\begin{equation}
\label{eqpsiPrime} \psi'(t) = \int V \,d(\nu-\mu_\star) + t
J_\alpha(\nu- \mu_\star) + J_\alpha(
\mu_\star,\nu- \mu_\star).
\end{equation}
Since $\mu_\star$ minimizes $I$, $\psi'(0^+)$ must be nonnegative:
\begin{eqnarray*}
0 &\leq&\int V \,d(\nu-\mu_\star) + J_\alpha(\mu_\star,
\nu-\mu_\star)
\\
&\leq&\int V \,d\nu+ J_\alpha(\mu_\star,\nu) - {{ \biggl(\int V
\,d
\mu_\star+ J_\alpha(\mu_\star) \biggr)}}
\\
&\leq&\int\bigl(V+U^{\mu_\star}_\alpha\bigr) \,d\nu-
C_\star.
\end{eqnarray*}
Therefore,
%
%
\begin{equation}
\label{eqconditionMinimalite} \forall\nu\in\mathcal{M}_1\cap
\mathcal{E}_\alpha
\qquad\int\bigl(V + U^{\mu_\star}_\alpha- C_\star\bigr) \,d
\nu\geq0.
\end{equation}
Since this holds for all $\nu$, $V + U^{\mu_\star}_\alpha$ is
greater than $C_\star$ quasi-everywhere.
Indeed, let $A = \{ x, V(x) + U^{\mu_\star}_\alpha(x) < C_\star\}$.
Since $V + U^{\mu_\star}_\alpha$ is measurable this is a Borel set.
Suppose by contradiction that its capacity is strictly positive.
By Proposition~\ref{propqe}, there exist a compact set $K\subset A$
and a measure $\nu$ with finite energy
supported in $K$. For this measure $\int V + U^{\mu_\star}_\alpha
\,d\nu
< C_\star$,
which contradicts (\ref{eqconditionMinimalite}).
This proves~(\ref{eqrob1}).

Let us prove (\ref{eqrob2}). Suppose $V(x) + U^{\mu_\star}_\alpha
(x) >
C_\star$
for some $x\in\operatorname{supp}(\mu_\star)$.
Since $V + U^{\mu_\star}_\alpha$ is lower semi-continuous,\vadjust{\goodbreak} we can find
a neighborhood $\mathcal{U}$ of $x$, and an $\eta>0$ such that
\[
\forall x \in\mathcal{U}\qquad V(x) + U^{\mu_\star}_\alpha(x)
\geq
C_\star+ \eta.
\]
Therefore
\[
\int\bigl(V + U^{\mu_\star}_\alpha\bigr) \,d\mu_\star
\geq(C_\star+ \eta) \mu_\star(\mathcal{U}) %
+ \int
_{\mathbb{R}^d\setminus\mathcal{U}} \bigl(V+ U^{\mu_\star}_\alpha
\bigr
) \,d
\mu_\star.
\]
Since $V+ U^{\mu_\star}_\alpha\geq C_\star$ quasi-everywhere, and
$\mu
_\star$
has finite energy, this holds $\mu_\star$ almost surely, so
\[
C_\star= \int V + U^{\mu_\star}_\alpha \,d\mu_\star
\geq C_\star+ \eta\mu_\star(\mathcal{U}).
\]
This is impossible since $\mu_\star(\mathcal{U})>0$, by definition
of the support.
Therefore (\ref{eqrob2}) holds.
\end{pf*}

\begin{pf*}{Proof of item (6) of Theorem~\ref{thriesz}}
Let
$\mu\in\mathcal{E}_\alpha\cap\mathcal{M}_1(\mathbb{R}^d)$ be
such that
$V\in\rL^1(\nu)$.
It is enough to show that, if (\ref{eqrob3prime}) and (\ref
{eqrob4prime}) hold,
then $\mu= \mu_\star$. We argue by contradiction and suppose $\mu
\neq
\mu_\star$.
Consider again the function $\psi$ (with $\nu=\mu$):
$\psi(t)=I((1-t)\mu_\star+ t \mu)$, $t\in[0,1]$.
According to Lemma~\ref{leconvexity}, this function is strictly convex,
therefore $\psi'(1) > \psi'(0) \geq0$. The explicit expression of
$\psi'$
[equation~(\ref{eqpsiPrime})] gives
\begin{eqnarray*}
0<\psi'(1) &=&\int V \,d(\mu- \mu_\star) +
J_\alpha(\mu- \mu_\star) + J_\alpha(\mu_\star,\mu- \mu_\star)
\\
&=& \int V \,d\mu- \int V \,d\mu_\star+ J_\alpha(\mu) -
J_\alpha(\mu,\mu_\star).
\end{eqnarray*}
Therefore,
%
%
\begin{equation}
\label{eqdracofeu} \int\bigl(U_\alpha^\mu+ V \bigr) \,d
\mu_\star< \int\bigl(U_\alpha^{\mu} + V \bigr) \,d\mu.
\end{equation}
On the other hand, integrating (\ref{eqrob3prime}) with respect to
$\mu$
and (\ref{eqrob4prime}) with respect to $\mu_\star$ yields
\[
\int\bigl(U_\alpha^\mu+ V \bigr) \,d\mu\leq C \leq\int
\bigl(U_\alpha^{\mu} + V \bigr) \,d\mu_\star,
\]
which contradicts (\ref{eqdracofeu}) and concludes the proof.
\end{pf*}

To prove the last result of Theorem~\ref{thriesz} we recall the following
classical result.

%
%
\begin{theorem}[(Principle of domination)]
Suppose $\alpha\leq2$. Let $\mu$ and $\nu$ be two positive measures in
$\mathcal{E}_\alpha$, and $c$ a nonnegative constant. If the inequality
\[
U_\alpha^\mu(x) \leq U_\alpha^\nu(x) + c
\]
holds $\mu$-almost surely, then it holds for all $x\in\mathbb{R}^d$.
\label{thmdomination}
\end{theorem}

\begin{pf}
In the Coulomb case $\alpha=2$, \cite{MR0350027}, Theorem~1.27, page~110,
applies, since $U_\alpha^\nu$ is positive and super-harmonic. If
$\alpha<2$,
the potential $U_\alpha^\nu$ is \mbox{$\alpha$-}superharmonic, so we can apply
\cite{MR0350027}, Theorem~1.29, page~115 and get the result.
\end{pf}

\begin{pf*}{Proof of item (7) of Theorem~\ref{thriesz}}
We follow the proof of Theorem~1.3 in~\cite{MR2276529}.
Arguing by contradiction, let us suppose that, for some measure $\mu$,
and some $\epsilon>0$,
\[
\sup_{\operatorname{supp}(\mu)} {{ \bigl(U_\alpha^\mu+ V
\bigr)}} \leq C_\star- \epsilon.
\]
By (\ref{eqrob2}), this implies that
\[
U_\alpha^\mu(x) + \epsilon\leq U_\alpha^{\mu_\star}(x),
\]
for all $x\in\operatorname{supp}(\mu)$. Let $\eta$ be the
equilibrium (probability)
measure of\break  $\operatorname{supp}(\mu)\dvtx U_\alpha^\eta(x) =
C_\eta$
on $\operatorname{supp}(\mu)$,
therefore
\[
U_\alpha^{\mu+ (\epsilon/C_\eta)\eta} \leq U_\alpha^{\mu_\star}
\]
for all $x$ in $\operatorname{supp}(\mu)$. By the principle of
domination this holds at
infinity. Since for any compactly supported $\mu$, $U_\alpha^\mu(x)
\sim
\frac{\mu(\mathbb{R}^d)}{{{\vert x \vert}}^{d-\alpha}}$ at infinity,
we get a contradiction
\[
(1+\epsilon/C_\eta) \leq1.
\]
Similarly, if
\[
\mbox{``}\inf_{\operatorname{supp}(\mu_\star)}\mbox{''}  \bigl(U_\alpha^\mu(x) +
V(x) \bigr) > C_\star,
\]
then $U_\alpha^\mu+ V \geq C_\star+ \epsilon$ q.e. on
$\operatorname{supp}(\mu_\star
)$, so
\[
U_\alpha^\mu(x) \geq U_\alpha^{\mu_\star}(x) +
\epsilon, \qquad\mu_\star\mbox{-a.s.}
\]
The same proof as before applies to get a contradiction.
\end{pf*}

\subsection{\texorpdfstring{Radial external fields in the Coulomb case: Corollary~\protect\ref{cocoulomb}.}{Radial external fields in the Coulomb case: Corollary~1.3}}\label{ssradialcoulomb}

For the sake of completeness, let us finally give a proof of the result
mentioned in Corollary~\ref{cocoulomb}.

Changing $V$ into $\beta V$, we can assume without loss of generality that
$\beta=1$.

Recall that $V$ is supposed to be radially symmetric and of class
$\mathcal{C}^2$: there exists $v\dvtx\mathbb{R}_+\to\mathbb{R}$
such that $V(x) =
v({{\vert x \vert}})$.

In this case it is thus natural to look for a radially symmetric equilibrium
probability measure. Guided by the results of \cite{MR1485778}, let us
consider an absolutely continuous probability measure $\mu$, such that
$\operatorname{supp}(\mu) = \{x\in\R^d; r_0\leq|x|\leq R_0\}$ for some
$0\leq
r_0<R_0$ and such that $d\mu= M(r) \,d\sigma_r \,dr$, where $M\dvtx[r_0,R_0]
\to
\mathbb{R}_+$ is assumed to be continuous.

First let us calculate the potential of $\mu$. Using the Gauss's averaging
principle (Theorem~\ref{thmgauss}), it holds
%
%
\begin{eqnarray}\label{eq=potentielDeMu}
U_2^\mu(x) &=& \int\!\!\!\int M(r) W(x,y) \,d
\sigma_r(y) \,dr\nonumber
\\
&=& \int M(r) \int_{\partial B(0,r)} \frac{1}{{{\vert x-y \vert
}}^{d-2}} \,d
\sigma_r(y)\,dr
\nonumber\\[-8pt]\\[-8pt]
&=& \sigma_d \int M(r) r^{d-1}{{ \biggl(
\frac{\IND_{{ \{{{\vert x \vert}}>r \}}}}{{{\vert x \vert
}}^{d-2}} + \frac
{\IND_{{ \{{{\vert x \vert}}\leq r \}
}}}{r^{d-2}} \biggr)}} \,dr\nonumber
\\
&=& \frac{\sigma_d}{{{\vert x \vert
}}^{d-2}} \int_0^{{\vert x
\vert}}
M(r) r^{d-1} \,dr + \sigma_d \int_{{{\vert x \vert}}}^\infty
M(r) r \,dr.\nonumber
\end{eqnarray}
Thus $U^\mu_2(x)=u(|x|)$, for some function $u$ of class $\mathcal{C}^1$.

Now, let us consider condition (\ref{eqrob3}). It holds if and only
if there
exists some $C$ such that $u(r)=C-v(r)$ for all $r\in[r_0,R_0]$. This is\vspace*{2pt}
obviously equivalent to the conditions $u'(r)=-v'(r)$ for all $r\in[r_0,R_0]$
and $u(R_0)=1/R_0^{d-2}=C-v(R_0)$ [here we use that $\sigma_d \int_0^{R_0}
M(r)r^{d-1} \,dr=1$]. Observing that
\[
u'(r)=-\frac{\sigma_d(d-2)}{r^{d-1}}\int_0^r
M(t)t^{d-1} \,dt,
\]
we see that $u'=-v'$ on $[r_0,R_0]$ if and only if $u'(R_0)=-v'(R_0)$ which
amounts to $w(R_0)=d-2$ and
$M(t)=\frac{1}{\sigma_d(d-2)}\frac{\omega'(t)}{t^{d-1}}$, for all
$t\in
[r_0,R_0]$, where we recall that $w(t)=t^{d-1}v'(t)$, $t\geq0$. The condition
$\sigma_d \int_0^{R_0} M(r)r^{d-1} \,dr=1$ implies that
$\frac{1}{d-2}(w(R_0)-w(r_0))=1$ and so $w(r_0)=0$. In the case where
$w$ is
increasing, this determines uniquely $r_0=0$ and $R_0=w^{-1}(d-2)$. In
the case
where $v$ is supposed to be convex, we see that $w$ is increasing on
$[a_0,\infty[$ with $a_0=\inf\{t >0; v'(t)>0\}$ and $w\leq0$ on $[0,a_0]$.
Therefore $R_0$ is uniquely defined and reasoning on the support of
$\mu$
easily yields to the conclusion that $r_0=a_0$. In all cases, the probability
$\mu$ is uniquely determined and $C=1/R_0^{d-2} +v(R_0)$.

It remains to check that this probability $\mu$ satisfies also condition
(\ref{eqrob4}). If $r=|x|\geq R_0$, then
$U_2^\mu(x)+V(x)=\frac{1}{r^{d-2}}+v(r) \geq\frac{1}{R_0^{d-2}}+v(R_0)=C$,
since it is easy to check that $r\mapsto\frac{1}{r^{d-2}} + v(r)$ is
increasing on $[R_0,\infty)$. In the case where $v$ is convex and
$r\leq r_0$,
an integration by parts yields to
\begin{eqnarray*}
U_2^\mu(x)&=&\frac{1}{d-2}\int_{r_0}^{R_0}
\frac{w'(t)}{t^{d-2}} \,dt=\frac{1}{d-2}\int_{r_0}^{R_0}
(d-1)v'(t)+tv''(t) \,dt
\\
&=&v(R_0)-v(r_0)+\frac{1}{R_0^{d-2}}=
C-v(r_0) \geq C-v(r),
\end{eqnarray*}
since $v$ is nonincreasing on $[0,r_0]$. Therefore, in all cases
$U^\mu_2(x)+V(x)\geq C$ for every $x\in\R^d$, which completes the
proof of
the characterization of the equilibrium measure.

Finally, if the external field $V$ is quadratic, that is, if $v(r) =
r^2$, then
$w(r) = 2r^d$ and so $r_0 = 0$, $R_0= ((d-2)/2)^{1/d}$ and $M(r) =
\frac{2\,d}{\sigma_d(d-2)}\IND_{{ \{{{\vert x \vert}}\leq
R_0 \}}}$. In other words, the
equilibrium probability measure is uniform on the ball centered in $0$
and of
radius $((d-2)/2)^{1/d}$.

\subsection{\texorpdfstring{Prescribed equilibrium measure.}{Prescribed equilibrium measure}}\label{ssproofcoprescription}

In this section we prove Corollary~\ref{coprescription}. We will need
the following elementary lemma.
%
%
\begin{lemma}[(Regularity of Riesz potential)]\label{lemprescription}
Let $0<\alpha<d$, $d\geq1$, and let $\mu$ be a probability measure
with a
density $f\in\mathbf{L}^p_{\mathrm{loc}}(\R^d)$ for some
$p>d/\alpha$. Then
$U_\alpha^{\mu}$ is continuous and finite everywhere on $\R^d$.
\end{lemma}

\begin{pf}
For all $n\geq1$, define
\[
R_n(x):=\int f(y) \min\bigl(n; k_{\Delta_\alpha}(x-y) \bigr) \,dy
\]
and
\[
S_n(x):=U_\alpha^\mu(x)-R_n(x)=\int
f(y) \bigl[k_{\Delta_\alpha
}(x-y)-n \bigr]_+ \,dy.
\]
It follows from the dominated convergence theorem that $R_n$ is continuous
on $\R^d$. Let us show that $S_n$ converges to $0$ uniformly on compact
sets, which will prove the claim. Let $q:=p/(p-1)$ be the conjugate exponent
of $p$; applying H\"older inequality yields to
\begin{eqnarray*}
0\leq S_n(x) &\leq&\int f(y)\frac{1}{|x-y|^{d-\alpha}}%
\IND_{B(x,n^{-1/(d-\alpha)})}(y) \,dy
\\
&\leq&\|f\|_{p,B(x,1)} %
{{ \biggl(\int\frac{1}{|x-y|^{q(d-\alpha)}}\IND
_{B(x,n^{-1/(d-\alpha)})}(y) \,dy \biggr)}}^{1/q}
\\
&=& \|f\|_{p,B(x,1)} \varepsilon_n,
\end{eqnarray*}
where\vspace*{2pt} $\varepsilon_n:=\sigma_d^{1/q} (\int_0^{n^{-1/(d-\alpha)}}
\frac{1}{u^{q(d-\alpha)-d+1}} \,du )^{1/q}$ and where $\sigma
_d$ is
the surface of the unit Euclidean ball. The condition $p>d/\alpha$ is
equivalent to $q(d-\alpha)-d+1<1$ and so $\varepsilon_n$ is finite
for all
$n$ and
$\varepsilon_n\to0$ as $n\to\infty$. We conclude from this that if
$K$ is a
compact set of $\R^d$ and $K_1=\{x\in\R^d; d(x,K)\leq1\}$, it holds
\[
\sup_{x\in K}|S_n|(x) \leq\|f\|_{p, K_1}
\varepsilon_n,
\]
which completes the proof.
\end{pf}

\begin{pf*}{Proof of Corollary~\ref{coprescription}}
Lemma~\ref{lemprescription} above shows that $U_\alpha^{\mu_\star
}$ is
continuous and everywhere finite on $\R^d$. Since $\mu_\star$ is
compactly supported, $U_\alpha^{\mu_\star}(x)\to0$ as $|x|\to
\infty$.
Therefore $V(x)\to\infty$, when $|x|\to\infty$. This proves \textup
{(H2)}. The other assumptions are straightforward.
By the very definition of $V$, it holds
\[
U_\alpha^{\mu_\star}(x)+V(x)\geq0\qquad\forall x\in\R^d,
\]
with equality on $B(0,R)\supseteq\operatorname{supp}(\mu)$.
According to
point (6) of Theorem~\ref{thriesz}, this proves that $\mu_*$ is the
(unique) minimizer of $I$. The last assertion follows from point~(4) of
Theorem~\ref{thriesz}.
\end{pf*}

\section*{\texorpdfstring{Acknowledgments.}{Acknowledgments}}
The authors would like to thank Luc Del\'{e}aval, Abey L\'{o}pez Garc\'
{i}a, Arnaud Guillin, Adrien Hardy, St\'ephane Mischler, and Karl-Theodor Sturm for
providing references and for interesting discussions, and also two anonymous
referees for their constructive remarks.


%

\printaddresses

\end{document}